# Uniqueness and non-uniqueness in percolation theory*

**Olle Häggström and Johan Jonasson**

*Dept of Mathematical Sciences*
*Chalmers University of Technology*
*412 96 Göteborg*
*SWEDEN*
*e-mail:* olleh@math.chalmers.se; jonasson@math.chalmers.se

**Abstract:** This paper is an up-to-date introduction to the problem of uniqueness versus non-uniqueness of infinite clusters for percolation on $\mathbb{Z}^d$ and, more generally, on transitive graphs. For iid percolation on $\mathbb{Z}^d$, uniqueness of the infinite cluster is a classical result, while on certain other transitive graphs uniqueness may fail. Key properties of the graphs in this context turn out to be amenability and nonamenability. The same problem is considered for certain dependent percolation models – most prominently the Fortuin–Kasteleyn random-cluster model – and in situations where the standard connectivity notion is replaced by entanglement or rigidity. So-called simultaneous uniqueness in couplings of percolation processes is also considered. Some of the main results are proved in detail, while for others the proofs are merely sketched, and for yet others they are omitted. Several open problems are discussed.

**AMS 2000 subject classifications:** Primary 60K35, 82B43.
**Keywords and phrases:** percolation, uniqueness of the infinite cluster, transitive graphs, amenability.

## 1. Introduction

Percolation theory began in 1957, when Broadbent and Hammersley [18] introduced the, now standard, iid bond percolation model on $\mathbb{Z}^d$. They did this to model a porous stone on a microscopic level in order to study the question of whether the center of the stone gets wetted when it is immersed into a bucket of water. Since then, percolation theory has attracted an enormous amount of interest. This is partly because of its applicability: the independence in the model makes it possible to handle, and still it is not too unrealistic for many random media. It has also proved fruitful to use comparisons with iid percolation in order to obtain results for various types of dependent percolation models. The main reason, however, for the interest is the mathematical beauty of the topic with its abundance of easily formulated conjectures of which many have proved to be very difficult to settle and which have inspired the development of powerful mathematical techniques.

*Research supported by the Swedish Research Council.

With the *comme il faut* slip of notation, we write $\mathbb{Z}^d$ for the graph whose vertex set is $\mathbb{Z}^d$ and whose edge set consists of all pairs of vertices sitting at Euclidean distance 1 from each other. While percolation theory began on $\mathbb{Z}^d$, the model makes perfect sense on any connected graph $G = (V, E)$: Let each edge be retained with probability $p$ (and deleted with the remaining probability $1-p$) independently of all other edges. On $\mathbb{Z}^d$ in $d \geq 2$ dimensions, it is a fundamental fact, dating back to Broadbent and Hammersley, that the occurrence of an infinite connected component – or an *infinite cluster*, as we will also call it – has probability 0 or 1 depending on whether $p$ is below or above a certain critical threshold $p_c \in (0, 1)$. When there are infinite clusters one may ask how many, and this is the topic of the present survey paper.

Let us briefly mention a motivational example from real life. Take $G = (V, E)$ be the graph whose vertex set consists of all mathematicians, and connect any two of them by an edge $e \in E$ if they have ever coauthored a mathematical paper. The distance in this graph between a given mathematician $v$ and Paul Erdős is colloquially referred to as $v$'s *Erdős number*. If we define a mathematician to be a person who has authored at least one piece of work that by May 2004 had found its way into the *MathSciNet* database, and two such mathematicians to be coauthors if they have a joint publication in that database, then the graph consists of about 401 000 vertices, and a good deal of other interesting statistics and graph characteristics are available; see Grossman [43, 44]. About 83 000 of the nodes are isolated, corresponding to mathematicians with no coauthors. The largest connected component – the one containing Erdős – contains about 268 000 vertices; these are all the mathematicians that have finite Erdős numbers. The *second* largest connected component contains only 32 vertices. This striking contrast between the size of the largest and the second largest connected component may be phrased as *uniqueness of the giant component*, and calls for our attention. Is this a special feature of social networks (as opposed to graph structures arising in other applications) or even of the social structure of mathematicians, or is it an instance of a much more general phenomenon?

Something similar happens in the well-known Erdős–Rényi random graph, consisting of $n$ vertices where each pair is linked by an edge with probability $p$, independently of all other pairs. (Equivalently, the Erdős–Rényi random graph can be seen as bond percolation with retention parameter $p$ on the complete graph with $n$ vertices.) There is an enormous literature on this model; see, e.g., Bollobás [15] or Janson et al. [70]. The most natural way to scale $p$ as $n \to \infty$ is to fix $c$ and let $p = \frac{c}{n}$, and the classical result going back to Erdős and Rényi [30] is that the proportion of vertices sitting in the largest connected component tends in probability to a constant, which is 0 when $c \leq 1$ and strictly positive when $c > 1$. On the other hand, it is fairly straightforward to show that for *all* $c$, the proportion of vertices sitting in the *second* largest component, tends in probability to 0. Hence, for $c > 1$ and $n$ large, we obtain another instance of the uniqueness of the giant component phenomenon.

In this paper, as in most of mathematical percolation theory, we will only be dealing with *infinite* graphs[1], where the issue of uniqueness of the giant component translates naturally into the question of whether there is a unique *infinite cluster*. This has the advantage of always having a clear-cut yes/no-answer, in contrast to the finite setting where it is not always totally obvious what one really should mean by a giant component.

From now on, $G = (V, E)$ will denote the underlying infinite graph on which the percolation process takes place. $G$ will always be assumed to be connected. Some more terminology: An iid *bond percolation* with retention parameter $p$ on the graph $G = (V, E)$ is a random element $X$ of $\{0, 1\}^E$ whose distribution $P_p$ is product measure with marginals $(1 - p, p)$. We identify $X$ with the subgraph of $G$ containing all vertices $v \in V$ and precisely those edges $e \in E$ taking value $X(e) = 1$. When $X(e) = 1$ we speak of $e$ as an *open* edge, whereas if $X(e) = 0$ we say that $e$ is *closed*.

In the equally natural process of iid *site percolation*, it is the vertices, rather than the edges, that are retained at random (independently, each with probability $p$). As far as qualitative results are concerned, it is usually of little importance whether bond or site percolation is considered. Here we will, with few exceptions, focus on bond percolation; most results and proofs have obvious analogues for site percolation.

The most basic question to ask about $X$ is whether it contains an infinite cluster. Write $C$ for this event, i.e., the event that $X$ contains at least one infinite cluster. It follows from the Kolmogorov 0-1-law that $P_p(C)$ is either 0 or 1. Furthermore, it is intuitively obvious that $P_p(C)$ should be non-decreasing in $p$, since adding edges to $X$ cannot destroy an infinite cluster. This intuition is turned into mathematical rigor via the following coupling, which will be useful to us several times later in this paper:

Coupling 1.1 The simultaneous coupling construction. *Associate with each edge $e \in E$ a random variable $U_e$ uniformly distributed on the unit interval. For each $p \in [0, 1]$, define $X_p(e) = I_{[0,p)}(U_e)$. Then, for any $p$, $X_p = \{X_p(e)\}_{e \in E}$ is an iid bond percolation with retention parameter $p$, i.e., $X_p$ has distribution $P_p$.*

For $0 \le p_1 < p_2 \le 1$, we have in the simultaneous coupling construction that $X_{p_1}(e) \le X_{p_2}(e)$ for every $e \in E$, and consequently $P_p(C)$ is non-decreasing in $p$. Combined with the observation that $P_p(C)$ must be 0 or 1, this implies the existence of a critical probability $p_c = p_c(G) \in [0, 1]$ such that

$$P_p(C) = \begin{cases} 0, & p < p_c \\ 1, & p > p_c \,. \end{cases}$$

For a given vertex $v \in V$, define $\theta_v(p) = P_p(v \leftrightarrow \infty)$, where $\{v \leftrightarrow \infty\}$ is the event that $v$ is in an infinite cluster. In many cases, such as when $G = \mathbb{Z}^d$ (and more generally when $G$ is transitive; see Definition 1.2 below), $\theta_v(p)$ is independent of the choice of $v \in V$, and we write $\theta_G(p)$ instead (or just $\theta(p)$ if it is obvious which graph is meant). Clearly, $P_p(C) = 0$ implies that $\theta_v(p) = 0$ for any $v$. Conversely, under our universal assumption that $G$ is connected, $P_p(C) > 0$ implies $\theta_v(p) > 0$ for any $v$; this is easy to see, for instance using Coupling 2.5 below. Hence, the following alternative characterization of the critical value $p_c(G)$ holds: for any $v \in V$, we have

$$\begin{cases} \theta_v(p) = 0, & p < p_c \\ \theta_v(p) > 0, & p > p_c\,. \end{cases}$$

[1]However, in recent years, there have been some very interesting developments in the project of bridging the infinite and the large-but-finite in percolation theory; see, e.g., [17] and [16].

When $p < p_c$ the percolation said to be *subcritical* while for $p > p_c$ it is said to be *supercritical*, and for $p = p_c$ it is said to be *critical*. Whether $P_p(C)$ is 0 or 1 at criticality depends on the choice of $G$. For $G = \mathbb{Z}^d$ with $d \geq 2$, it is believed that $\theta_G(p_c) = 0$. This is known only for $d = 2$, where it was established by the work of Harris [62] and Kesten [76], and for $d \geq 19$ due to Hara and Slade [61]. The cases $3 \leq d \leq 18$ are what remain of the most classical long-standing open problem in percolation theory. Benjamini and Schramm [12] (we will soon hear more about their paper) extended the conjecture $\theta_G(p_c) = 0$ to the class of transitive graphs with $p_c < 1$; this was later shown to be the case for a large class of graphs by Benjamini et al. [9] – see Theorem 7.3 below – but the full conjecture remains open.

We move on to the main topic of this paper: in supercritical percolation, how many infinite clusters are there? A very brief history of this problem is as follows.

In his classical 1960 paper, Harris [62] showed that on $\mathbb{Z}^2$ one has $p_c \geq 1/2$, and by combining this with the self-duality of $\mathbb{Z}^2$ he was able to deduce that for any supercritical $p$ there is a.s. a unique infinite cluster. The world then had to wait until 1987 before Aizenman et al. [3] were able to extend this to $d \geq 3$ and thus establish uniqueness of the infinite cluster in all dimensions. A couple of years later, Burton and Keane [19] obtained an alternative proof which is much shorter and easier to comprehend, and more amenable to generalizations. Consequently, the Burton–Keane proof has become the standard one to present in textbooks and courses (see, e.g., Grimmett [37]), and we will follow suit.

Soon afterwards, the subject was given an additional spark in a paper by Grimmett and Newman [42] who showed that there are interesting graphs for which one has uniqueness of the infinite cluster for some $p$ but not for others. The specific example they used was the Cartesian product graph $G = \mathbb{T}_d \times \mathbb{Z}$ where $\mathbb{T}_d$ is a regular tree of degree $d + 1$, and they showed that the desired phenomenon occurs when $d$ is large enough: when $p$ is sufficiently close to $p_c$, there are infinitely many infinite clusters, while uniqueness holds when $p$ is sufficiently close to 1. The Grimmett–Newman paper can be seen as a forerunner to the highly influential paper from 1996 by Benjamini and Schramm [12], who suggested – correctly, as it turned out – that a fruitful generality in which to study percolation would be the *transitive* and the *quasi-transitive* graphs; see the following definition. A frenzy of activity and interesting results followed upon the Benjamini–Schramm paper, including a number of results concerning uniqueness versus non-uniqueness of infinite clusters.

Definition 1.2 *Let $G = (V, E)$ be an infinite graph. A bijective map $f : V \to V$ such that $\{f(u), f(v)\} \in E$ if and only if $\{u, v\} \in E$ is called a* **graph automorphism** *for $G$. The graph $G$ is said to be* **transitive** *if for any $u, v \in V$ there exists a graph automorphism mapping $u$ on $v$. More generally, $G$ is said to be* **quasi-transitive** *if $V$ can be partitioned into a finite number of vertex sets $V_1, \dots V_k$ such that for any $i \in \{1, \dots, k\}$ and any $u, v \in V_i$, there exists a graph automorphism mapping $u$ on $v$.*

Heuristically, the graph $G$ is transitive if and only if it "looks the same" as seen from any vertex – examples include both the usual $\mathbb{Z}^d$ lattice and the Grimmett–Newman example $\mathbb{T}_d \times \mathbb{Z}$ – while quasi-transitivity means that there are only finitely many "kinds" of vertices.

The topic of uniqueness (and non-uniqueness) of infinite clusters in percolation theory was reviewed already in a 1994 paper entitled *Uniqueness in percolation theory* by Meester [83], but the subject has, as we shall see, developed vastly since then. The title of the present paper is meant as an allusion to Meester's paper, and as an indication of the shift in the subject's center of mass that has taken place during the last decade. See also Grimmett [39] for a more recent treatment, with a somewhat different emphasis compared to ours.

The rest of this paper is organized as follows. In Section 2, we consider the $\mathbb{Z}^d$ case, and give the uniqueness arguments of Harris [62] and Burton and Keane [19]. In Section 3, we outline the more recent work concerning percolation on more general transitive (and quasi-transitive) graphs following the footsteps of Grimmett and Newman [42] and Benjamini and Schramm [12]. Sections 4–9 are then devoted to describing this work in more detail, including a number of proofs. In particular, we will learn about the crucial role played by the properties of *amenability* versus *nonamenability* of transitive graphs, and the progress towards the still-open problem (see Conjecture 3.3) of establishing whether nonamenability is equivalent to the existence of some $p \in (0, 1)$ where non-uniqueness holds.

Then, in Section 10, we move on to questions concerning so-called *simultaneous uniqueness* in coupling constructions such as the one in Coupling 1.1, and in Section 11 we briefly mention similar issues for so-called *dynamical percolation.* In Sections 12 and 13, we consider what happens when we let go of the iid assumption in favor of dependent percolation models; prime examples are the random-cluster model and uniform spanning trees. Some particular results for dependent percolation on $\mathbb{Z}^2$ are considered in Section 14. Finally, in Sections 15 and 16, we consider what happens if we focus not on connected components, but rather on *entangled* or on *rigid* components.

Before moving on, let us finally emphasize that we make no claims at a complete coverage of the topic of uniqueness vs non-uniqueness of infinite clusters in percolation theory – the subject is much too large for that. In our choice of material, we have (apart from some early fundamentals of the subject) mainly been guided by our wish to include those topics that we feel have developed most vigorously during the last decade, but the choices are of course open to criticism. Two topics that were treated at some length by Meester [83] but are completely omitted here, are continuum percolation [84] and fractal percolation [24].

## 2. Percolation on $\mathbb{Z}^d$

We begin this section with Harris' [62] proof of uniqueness of the infinite cluster for iid bond percolation on $\mathbb{Z}^2$. For this, we will need to consider the *dual* of $\mathbb{Z}^2$. For later use, we take the opportunity to define the dual graph of a planar graph in general.

Definition 2.1 *Let $G = (V, E)$ be a planar graph with a specified embedding in $\mathbb{R}^2$. The dual graph $G^\dagger = (V^\dagger, E^\dagger)$ is the graph in which $V^\dagger$ is the set of faces of $G$ and for $u, v \in V^\dagger$, $\{u, v\} \in E^\dagger$ iff the faces $u$ and $v$ share an edge in $E$.*

For $G = \mathbb{Z}^2$, it is natural to picture the dual as a copy of the original $\mathbb{Z}^2$-lattice, shifted in $\mathbb{R}^2$ by the vector $(\frac{1}{2}, \frac{1}{2})$.

Now, if $X$ is an iid bond percolation with retention parameter $p$ on $G$, then the *dual percolation* $X^\dagger$ on $G^\dagger$ is defined by letting a given edge in $G^\dagger$ be open in $X^\dagger$ iff its crossing edge in $G$ is closed in $X$. By this definition, $X^\dagger$ is an iid bond percolation on $G^\dagger$ with retention parameter $1-p$. The point of considering dual percolation for $\mathbb{Z}^2$ is that the connected component of a vertex $v \in \mathbb{Z}^2$ is finite iff $v$ is surrounded by a *circuit*, i.e., a cycle of open edges, in $X^\dagger$. This statement is easy to believe but cumbersome to prove in full rigor (see, e.g., Kesten [77]); we state it here without proof.

Harris' result on uniqueness is the following.

Theorem 2.2 (Harris [62]) *The infinite cluster of supercritical percolation on $\mathbb{Z}^2$ is a.s. unique.*

*Proof.* By an ingenious sequence of geometric arguments that we do not present here, Harris showed that on $\mathbb{Z}^2$ percolation does not occur for $p = 1/2$.[2] The dual of $\mathbb{Z}^2$ is isomorphic to $\mathbb{Z}^2$ itself, and if one takes $p = 1/2$, the dual percolation $X^\dagger$ has exactly the same probabilistic behavior as the original percolation $X$. Let $B$ be a large box in $\mathbb{Z}^2$ centered at the origin. Since a.s. there is no infinite path of open edges in $X$ intersecting $B$, $B$ is a.s. surrounded by a circuit in $X^\dagger$. Thus $B$ must also be surrounded by a circuit in $X$. Using Coupling 1.1, we realize that this also holds for all higher values of $p$, in particular for all $p$ for which percolation occurs. This means that no two separate infinite clusters can both intersect $B$, because they would then have to be be joined by the circuit. Since the size of $B$ was arbitrary, this shows that there can only be one infinite cluster. □

Harris' argument, although beautiful, is obviously highly dependent on the planar structure of $\mathbb{Z}^2$, and it would be a waste of time to try to generalize the argument in order to prove, e.g., uniqueness of the infinite cluster for percolation on $\mathbb{Z}^d$ for $d \geq 3$. More than 25 years later, Aizenman et al. [3] were able to prove uniqueness of the infinite cluster for percolation on $\mathbb{Z}^d$ in arbitrary dimension $d$. Their proof was rather difficult, but a substantially simpler proof was soon obtained by Burton and Keane [19]. The Burton–Keane proof exploits only two aspects of the graph structure of $\mathbb{Z}^d$: transitivity (recall Definition 1.2) and so-called *amenability*. Amenability will continue to play an important role in later sections, so we will state and prove the uniqueness result in that generality.

[2] Consequently $p_c \geq 1/2$. It quickly became a famous open problem to show that this inequality for $p_c$ is in fact an equality. 20 years later, Kesten [76] proved that $p_c$ indeed equals 1/2.

To define amenability, we first need a notion of isoperimetric constants. For an infinite connected graph $G = (V, E)$, define its edge-isoperimetric constant as

$$\kappa_E(G) = \inf_W \frac{|\partial_E W|}{|W|}\,,$$

where the infimum ranges over all finite connected subsets $W$ of $V$, and $\partial_E W$ is the set of edges with one end-vertex in $W$ and one in $V \setminus W$.

Definition 2.3 *Take $G$ to be a graph with edge-isoperimetric constant $\kappa_E(G)$. If $\kappa_E(G) = 0$ then $G$ is said to be* **amenable**, *while if $\kappa_E(G) > 0$ then we say that $G$ is* **nonamenable**.

We remark that it is equally natural to define the (inner) vertex-isoperimetric constant as

$$\kappa_V(G) = \inf_W \frac{|\partial_V W|}{|W|}$$

where $\partial_V W$ is the set of vertices of $W$ with at least one neighbor in $V \setminus W$. If $G$ has bounded degree (in particular, if $G$ is transitive), then $\kappa_V(G) = 0$ iff $\kappa_E(G) = 0$, so that the notion of amenability is independent of whether it is defined using the edge-isoperimetric or the vertex-isoperimetric constant.

To see that the $\mathbb{Z}^d$ lattice is amenable, it suffices to take $W_n = \{-n, \ldots, n\}^d$ and note that $\frac{|\partial_E W_n|}{|W_n|}$ tends to 0 as $n \to \infty$. Thus, uniqueness of the infinite cluster for percolation on $\mathbb{Z}^d$ is a special case of the following result.

Theorem 2.4 (Burton and Keane [19]) *Assume that $G$ is connected, transitive and amenable. Then, for any $p \in (0, 1)$ such that $P_p(C) > 0$, we get $P_p$-a.s. a unique infinite cluster.*

Like many arguments in this field – we shall see a number of examples later on in this paper – the Burton–Keane proof makes repeated use of a technique colloquially known as *local modification*. Perhaps the clearest way to describe the technique is via the following coupling.

Coupling 2.5 The local modifier. *Fix $p \in (0, 1)$ and an infinite graph $G = (V, E)$. Let $\mathcal{E}$ denote the set of all finite subsets of $E$. The set $\mathcal{E}$ is countable, and it is therefore possible to find a probability measure on $\mathcal{E}$ which assigns positive probability to all elements; let $Q$ be such a probability measure. Two $\{0, 1\}^E$-valued random objects $X$ and $X'$, each having distribution $P_p$, may be obtained as follows.*

1. *Pick a finite edge set $F \in \mathcal{E}$ randomly with distribution $Q$.*
2. *Conditionally on Step 1, and for each edge $e \in E \setminus F$ independently, set $X(e) = X'(e) = 0$ or $1$ with respective probabilities $1-p$ and $p$.*
3. *Conditionally on Steps 1 and 2, and for each edge $e \in F$ independently, set $X(e) = 0$ or $1$, and independently $X'(e) = 0$ or $1$, each with probability $p$ of taking value $1$.*

This coupling was explicitly introduced by Häggström [52], but arguments based on variants of the coupling had appeared implicitly many times before. It may at first sight look innocuous, but we shall soon see examples of its usefulness. The randomization in Step 1 is, for many applications, not necessary, but we state the coupling with this randomization in order to maximize its one-size-fits-all quality.

An important preliminary step towards Theorem 2.4 is the following result of Newman and Schulman [86].

LEMMA 2.6 *Fix $p \in [0,1]$ and a transitive graph $G = (V, E)$. The number of infinite clusters arising in iid bond percolation on $G$ with retention parameter $p$ is an a.s. constant, taking one of the values $0$, $1$ and $\infty$.*

*Proof.* Define the random variable $N$ as the number of infinite clusters, and for $n \in \{0, 1, 2, \ldots\} \cup \{\infty\}$ let $D_n$ denote the event that $N = n$. We first show that $N$ is an a.s. constant. Here we may appeal to the shift-invariance of $D_n$ and the well-known fact that an iid process is ergodic – meaning that all shift-invariant events have probability 0 or 1. Readers that are happy with this may fast-forward to the paragraph containing eq. (2.4). Others are invited to the following more explicit argument, which in fact is good enough to prove ergodicity (simply by replacing $D_n$ by an arbitrary shift-invariant event).

Assume for contradiction that there is an $n$ such that

$$0 < P_p(D_n) < 1\,, \tag{2.1}$$

and fix such an $n$. For $W \subset V$, write $X(W)$ as shorthand for $\{X(v)\}_{v \in W}$. For $v \in V$ and a positive integer $k$, let $B(v,k)$ denote the set of vertices sitting within (graph-theoretic) distance $k$ from $v$, and define the analogous edge set

$$B_E(v,k) = \{e \in E : \text{both endpoints of } e \text{ are in } B(v,k)\}\,.$$

Define the random variable $I_{D_n}$ as the indicator of the event $D_n$. Furthermore, for any integer $k$ and $v \in V$, define $I_{n,v,k}$ as "the best guess" of $I_{D_n}$ given $X(B_E(v,k))$, by which we mean that

$$I_{n,v,k} = \begin{cases} 0 & \text{if } P_p(D_n | X(B_E(v,k))) \le 1/2 \\ 1 & \text{if } P_p(D_n | X(B_E(v,k))) > 1/2\,. \end{cases}$$

It is an immediate consequence of Lévy's 0-1-law (see, e.g., [27, Sect. 4.5]) that for fixed $v$,

$$\lim_{k \to \infty} I_{n,v,k} = I_{D_n} \text{ almost surely.} \tag{2.2}$$

Next let $w_1, w_2 \ldots$ be a sequence of vertices such that for each $k$, $w_k$ sits at distance at least $2k$ from $v$ (the point of this being that $B_E(v,k)$ and $B_E(w_k,k)$ do not intersect). Clearly, the pair $(I_{D_n}, I_{n,w_k,k})$ has the same joint distribution as $(I_{D_n}, I_{n,v,k})$, and we therefore get from (2.2) that $I_{n,w_k,k}$ converges in probability to $I_{D_n}$. Combining this with (2.2) again yields

$$\lim_{k\to\infty} P_p(I_{n,w_k,k} = I_{n,v,k} = I_{D_n}) = 1\,. \tag{2.3}$$

On the other hand, $I_{n,v,k}$ and $I_{n,w_k,k}$ are independent, because they are defined on disjoint edge sets, and it follows using the assumption (2.1) that

$$\lim_{k\to\infty} P_p(I_{n,w_k,k} = 1 = 1 - I_{n,v,k}) = P_p(D_n)(1 - P_d(D_n)) > 0$$

But this contradicts (2.3), so the assumption (2.1) must be false, and we have shown that the number of infinite clusters $N$ is an a.s. constant.

It remains to rule out that this a.s. constant equals some $n \in \{2, 3, \ldots\}$. Fix such an $n$, and also a vertex $v \in V$, and assume for contradiction that

$$P_p(N = n) = 1\,. \tag{2.4}$$

Then (since $G$ is connected) there exists a $k$ such that with positive probability $B(v,k)$ is intersected by all $n$ infinite clusters. Now pick $X, X' \in \{0,1\}^E$ and $F \in \mathcal{E}$ according to the local modifier (Coupling 2.5). By the choice of $k$, we have with positive probability that all infinite clusters in $X$ intersect $B(v,k)$. Conditional on that event, we have with positive probability that $F = B_E(v,k)$ and $X'(B_E(v,k) \equiv 1)$. But if these things happen, then $X'$ has a unique infinite cluster, so that

$$P_p(N = 1) > 0$$

contradicting (2.4), as desired. Hence, $N$ must be 0, 1 or $\infty$. $\Box$

*Proof of Theorem 2.4.* We will essentially follow [37, Sect. 8.2]; the proof there is on $\mathbb{Z}^d$ but is easily translated to the current setting. By Lemma 2.6 it suffices to rule out the possibility that $X$ contains infinitely many infinite clusters. First note that it may be assumed that the degree of the vertices in $G$ is at least 3 since otherwise $G$ would be isomorphic to $\mathbb{Z}$, a trivial case. We say that a vertex $v \in V$ is a *trifurcation* if

(a) $v$ is in an infinite cluster,
(b) $v$ is incident to exactly three open edges, and
(c) the removal of $v$ divides its infinite cluster into exactly three disjoint infinite cluster (and thereby no finite clusters).

Suppose that $A$ is a finite set of trifurcations belonging to the same infinite cluster $K$. Say that a member of $A$ is an *outer* member if at least two of the disjoint infinite clusters resulting from its removal contain no other member of $A$. We claim that $A$ must contain some outer member. To prove this, pick a member $v_1$ of $A$. If $v_1$ is not outer, then at least two of the disjoint infinite clusters resulting from its removal contain other members of $A$. Let $v_2$ and $v_3$ be such members. Now consider $v_3$. If $v_3$ is outer we are done. If not, the removal of $v_3$ results in three disjoint infinite clusters of which exactly one contains $v_1$ and $v_2$, and one contains some other member, $v_4$, of $A$. Now repeat for $v_4$: If $v_4$ is outer we are done. If not, its removal results in three disjoint infinite clusters of which exactly one contains $v_1, v_2$ and $v_3$, and one contains some other member, $v_5$, of $A$. Keep repeating this procedure until an outer member is found. This will happen sooner or later since $A$ is finite.

Next we claim that if $A$ is again a finite set of trifurcations in the same infinite cluster, then the removal of all of them will divide the cluster into at least $|A|+2$ disjoint infinite clusters. We do this by induction on $|A|$. The claim is obviously true when $|A| = 1$, so assume that it also holds for $|A| = j$ and consider a set $A$ of $j + 1$ trifurcations in the same infinite cluster. Let $v$ be an outer member of $A$. Remove all vertices of $A \setminus \{v\}$. This divides the infinite cluster into at least $j + 2$ disjoint ones. Now remove $v$. Since $v$ is outer, this results in an extra infinite cluster as desired.

Thus an infinite cluster with $j$ trifurcations in a finite set $W \subset V$ must intersect $\partial W$ in at least $j + 2$ vertices. Consequently, $W$ cannot contain more than $|\partial W| - 2$ trifurcations. Now let $t$ denote the probability that a vertex $v \in V$ is a trifurcation, and let $T(W)$ denote the number of trifurcations in $W$. Since $G$ is transitive, $t$ does not depend on the particular vertex. Therefore $\mathbb{E}T(W) = |W|t$. On the other hand, $T(W)$ is bounded by $|\partial W| - 2$, and so

$$t \leq \frac{|\partial W| - 2}{|W|}.$$

Since $G$ is amenable, the ratio on the right hand side can be made arbitrarily small by a suitable choice of $W$. Thus $t = 0$.

For the punch-line, we use the local modifier again. Pick $X, X' \in \{0,1\}^E$ and $C \in \mathcal{E}$ according to Coupling 2.5. Assume for contradiction that $N = \infty$ with $P_p$-probability 1. We can then find a finite connected vertex set $W$ such that with positive probability, at least three of the infinite clusters in $X$ intersect $W$. Conditional on this event, we have positive probability for the event that $F = E_W$, where $E_W$ is the set of edges $e \in E$ whose endpoints are both in $W$. Conditional on these events, each outcome of $X'(E_W)$ has positive probability, and it is easy to see that some of them consist of three disjoint simple paths, each emanating from one of these infinite clusters (a different one for each of the three paths), that are disjoint except that they all end at the same vertex $v \in V$. If $X'(E_W)$ happens to be such a configuration, then $v$ is a trifurcation for $X'$. Hence, $v$ has positive probability of being a trifurcation, so that $t > 0$, which is a contradiction. □

## 3. Percolation on more exotic graphs

In this section, we briefly outline the history of uniqueness versus non-uniqueness of infinite clusters for percolation beyond the classical $\mathbb{Z}^d$ setting, postponing proofs and other details to Sections 4–9.

Having come as far as establishing uniqueness of the infinite cluster on $\mathbb{Z}^d$ (and more generally on amenable transitive graphs), it is natural to ask how much further this result extends. Perhaps to all (connected) graphs? The answer is no, even if we restrict to transitive graphs, and the most straightforward way to see this is to consider the case when $G$ is a tree. Recall from Section 1 that $\mathbb{T}_d$ denotes the infinite tree of degree $d+1$, i.e., for the (unique) connected acyclic graph whose vertices all have degree $d+1$. The critical value $p_c(\mathbb{T}_d)$ for iid bond percolation on $\mathbb{T}_d$ satisfies $p_c(\mathbb{T}_d) = 1/d$, which is easy to show using, e.g., a branching process argument. It is then not a difficult task to deduce that for $d \geq 2$ and any $p \in (1/d, 1)$, the number of infinite clusters is a.s. $\infty$.

Thus, we have established the existence of two types of transitive graphs: those – such as $\mathbb{Z}^d$ with $d \geq 2$ – for which percolation with $p \in (p_c, 1)$ always yields a unique infinite cluster, and those – such as $\mathbb{T}_d$ with $d \geq 2$ – for which percolation with $p \in (p_c, 1)$ always yields infinitely many infinite clusters. Lemma 2.6 leaves us with one more possibility: perhaps there are graphs which yield a unique infinite cluster for *some* $p \in (p_c, 1)$, and infinitely many infinite clusters for *other* values of $p \in (p_c, 1)$?

Interestingly enough, it turns out that there do exist transitive graphs exhibiting this considerably more intricate behavior. The first example was found in a 1990 paper by Grimmett and Newman [42]. They considered the Cartesian product graph $G = \mathbb{T}_d \times \mathbb{Z}$, where Cartesian products are defined as follows.

Definition 3.1 *Let $G_1 = (V_1, E_1)$ and $G_2 = (V_2, E_2)$ be two graphs. The Cartesian product $G_1 \times G_2$ has vertex set $V_1 \times V_2$ and $\{(u_1, u_2), (v_1, v_2)\}$ is an edge of $G_1 \times G_2$ iff either $u_1 = v_1$ and $\{u_2, v_2\} \in E_2$ or $u_2 = v_2$ and $\{u_1, v_1\} \in E_1$.*

Thus a product graph $G_1 \times G_2$ has one "$G_1$-dimension" and one "$G_2$-dimension", and we refer to an edge $\{(u_1, u_2), (v_1, v_2)\}$ as a $G_1$-edge if $u_2 = v_2$ and a $G_2$-edge if $u_1 = v_1$. A basic example of a product graph is $\mathbb{Z} \times \mathbb{Z}$, which is identical to the square lattice $\mathbb{Z}^2$.

Grimmett and Newman proved for $d \geq 6$ that in iid percolation on $\mathbb{T}_d \times \mathbb{Z}$ in the supercritical phase $p \in (p_c, 1)$, the number $N$ of infinite clusters satisfies $N = \infty$ a.s. for $p$ sufficiently close to $p_c$, and $N = 1$ for $p$ sufficiently close to 1:

Theorem 3.2 (Grimmett and Newman [42]) *Let $d \geq 6$, let $G = \mathbb{T}_d \times \mathbb{Z}$, and consider an iid bond percolation $X$ on $G$ with retention parameter $p$. We then have $p_c < 1$ and the existence of $p_1$ and $p_2$ such that $p_c < p_1 \leq p_2 < 1$ such that*

- *for $p \in (p_c, p_1)$, $X$ contains a.s. infinitely many infinite clusters, whereas*
- *if $p \in (p_2, 1]$, $X$ contains a.s. a unique infinite cluster.*

The natural condition for this result is $d \geq 2$ rather than $d \geq 6$, and as later remarked by Schonmann [94], it can in fact be shown to hold in that greater generality using a combination of ideas from [94] and from Stacey [98]. The result is easier to prove for large $d$; in Section 4 we will outline a proof for $d \geq 65$.

Given Theorem 3.2 and Lemma 2.6, it is natural to conjecture that it should be possible in the statement of Theorem 3.2 to collapse $p_1$ and $p_2$ into a single value, uniquely separating the $N = \infty$ and the $N = 1$ parts of the supercritical regime. This would rule out the unlikely-sounding alternative that $N$ switches back and forth between $\infty$ and 1 as $p$ increases from $p_1$ to $p_2$. Letting

$$p_u = p_u(G) = \inf\{p : P_p(\text{there is a unique infinite cluster}) = 1\},$$

one would thus believe that when $p > p_u$, one has a.s. uniqueness of the infinite cluster. This is called *uniqueness monotonicity*, and was established (in much greater generality) in a pair of papers published back-to-back in the same journal issue in 1999: Häggström and Peres [56] proved uniqueness monotonicity for transitive (and quasi-transitive) graphs under the additional assumption of so-called *unimodularity* (see Definition 5.2 below), and Schonmann [95] gave a different proof where the unimodularity condition was shown to be superfluous. See Section 5.

Thus for all $G$ there exists a threshold $p_u \in [p_c, 1]$ such that uniqueness of the infinite cluster holds for all $p > p_u$, but for no $p < p_u$. At this point, a number of natural questions arise: For which graphs do we have $p_u > p_c$? For which graphs do we have $p_u < 1$? And what happens *at* $p_u$? (These questions and others were asked, in the setting of transitive and quasi-transitive graphs, in the 1996 paper by Benjamini and Schramm [12].)

For the regular tree $\mathbb{T}_d$ we saw that $p_u = 1$. The uniqueness results in Section 2 imply that $p_u = p_c$ for $G = \mathbb{Z}^d$ in any dimension $d$, and more generally when $G$ is transitive and amenable. The Grimmett–Newman example exhibits the more interesting behavior that $p_c < p_u < 1$. Benjamini and Schramm [12] conjectured the following:

Conjecture 3.3 *Let $G$ be any infinite connected transitive graph. Then $p_u > p_c$ if and only if $G$ is nonamenable.*

The Burton–Keane argument (i.e., the proof of Theorem 2.4) took care of the 'only if' part, so if the conjecture is true, then the Burton–Keane argument is in a rather deep sense sharp. If true, the conjecture would also provide a beautiful percolation characterization of amenability (analogous characterizations in terms of random walks and the Ising model have been obtained by, respectively, Kesten [74, 75] and Jonasson and Steif [72]). The 'if' part of the conjecture remains, in its full generality, open, but progress has been made in a number of directions: Lalley [78] proved that $p_u > p_c$ for certain planar Cayley graphs. This was improved by Benjamini and Schramm [13] who verified the conjecture for all planar graphs with the property of having one *end* (see the following definition). In a different direction, Pak and Smirnova-Nagnibeda [88] proved that $p_u > p_c$ for the Cayley graph of any nonamenable group provided the set of generators is properly chosen. These results are discussed in further detail in Sections 6 and 7. We consider proving (or disproving!) Conjecture 3.3 in its full generality to be one of the main open problems in this area.

Moving on to the problem of determining when $p_u < 1$, the following concept turns out to be of some importance.

DEFINITION 3.4 *Let $G = (V, E)$ be an infinite connected graph and for $W \subset V$ let $N_W$ be the number of infinite clusters the removal of $W$ from $G$ results in. The number $\sup_W N_W$, where the supremum is taken over all finite $W$, is called the* **number of ends** *of $G$.*

For example, $\mathbb{Z}$ has two ends, $\mathbb{Z}^d$ for $d \geq 2$ has one end, and the regular tree $\mathbb{T}_d$ with $d \geq 2$ has infinitely many ends. The Grimmett–Newman example $G = \mathbb{T}_d \times \mathbb{Z}$ has one end. Babson and Benjamini [6] showed that $p_u < 1$ holds on a certain class of one-ended graphs; see Section 8. This contrasts with the situation on a regular tree, where one has $p_u = 1$ (but also infinitely many ends).

Finally, there is the issue of what happens at the uniqueness critical value $p_u$; this question is most pertinent when $p_c < p_u < 1$. Do we get $N = 1$ or $N = \infty$? Somewhat surprisingly, the answer turns out to depend on $G$. Schonmann [94] proved that on $\mathbb{T}_d \times \mathbb{Z}$ the number of infinite clusters at $p_u$ is $\infty$, and Peres [90] generalized this by showing that for products $G \times H$ of graphs where at least one of $G$ and $H$ is nonamenable and both are infinite, the number of infinite clusters at $p_u$ cannot be 1. In contrast, Benjamini and Schramm [13] showed that for planar nonamenable graphs with one end, there is a unique infinite cluster at $p_u$. See Section 9 for further discussion of these results. To determine in general for which nonamenable transitive graphs uniqueness holds at $p_u$ is a highly interesting, but possibly quite difficult, open problem.

## 4. The Grimmett–Newman example $G = \mathbb{T}_d \times \mathbb{Z}$

In this section we outline a proof, under a stronger condition on $d$, of the Grimmett–Newman result that the supercritical regime of iid percolation on $\mathbb{T}_d \times \mathbb{Z}$ has both a uniqueness and a non-uniqueness part:

PROPOSITION 4.1 *Let $d \geq 17$, let $G = \mathbb{T}_d \times \mathbb{Z}$ and consider an iid bond percolation $X$ on $G$ with retention parameter $p$. If $p \in (1/d, 1/(4\sqrt{d+1}))$ then $X$ a.s. contains infinitely many infinite clusters, whereas if $p > 1/2$, then $X$ a.s. contains a unique infinite cluster.*

Note, to understand the relevance of the $d \geq 17$ condition, that $d = 17$ is the smallest integer for which $1/d < 1/(4\sqrt{d+1})$.

The proof requires a couple of lemmas, and for those we need some more notation and background. For percolation on a graph $G$ we denote for two vertices $u, v \in G$ by $\{u \leftrightarrow v\}$ the event that there is a path of open edges between $u$ and $v$. For two disjoint vertex sets $W, W' \subset V$, we write $W \leftrightarrow W'$ for the event that $\{u \leftrightarrow v\}$ occurs for some $u \in W$ and $v \in W'$.

An event $A \subset \{0,1\}^E$ is said to be *increasing* if for all $\xi, \xi' \in \{0,1\}^E$ such that $\xi(e) \le \xi'(e)$ for every $e \in E$ and $\xi \in A$, we also have $\xi' \in A$. Recall the Harris–FKG inequality (see, e.g., Harris [62] or Grimmett [37]), which says that under the product measure $P_p$, the indicator functions of any two increasing events are positively correlated.

LEMMA 4.2 *Suppose that for a given p, percolation on a transitive graph G produces a unique infinite cluster. Then there exists an $a > 0$ such that $P_p(u \leftrightarrow v) \ge a$ for every pair of vertices of G.*

*Proof.* Obviously, $\theta(p) > 0$. By the Harris–FKG inequality and the assumption of uniqueness of the infinite cluster,

$$\begin{aligned} P_p(u \leftrightarrow v) &\ge P_p(\{u \leftrightarrow \infty\} \cap \{v \leftrightarrow \infty\}) \\ &\ge P_p(u \leftrightarrow \infty) P_p(v \leftrightarrow \infty) \\ &\ge \theta(p)^2 , \end{aligned}$$

so putting $a = \theta(p)^2$ proves the lemma. □

LEMMA 4.3 (SCHONMANN [95]) *Suppose that iid percolation on a transitive graph G for a given p behaves in such a way that for every $a > 0$ there is an N sufficiently large so that $P_p(B(u,N) \leftrightarrow B(v,N)) \ge 1 - a$ for every u and v. Then there is a.s. a unique infinite cluster.*

This result is intuitively easy to believe: if, for fixed $u$, $P_p(B(u,N) \leftrightarrow B(v,N))$ is bounded away from 0, then with positive probability $B(u,N)$ is intersected by an infinite cluster which is sufficiently "omnipresent" that it appears to rule out the existence of other infinite clusters. Or to put it in another way: if the a.s. number of infinite clusters is 0 or $\infty$, then some small parts of the graph sitting sufficiently far from each other should have a very hard time connecting to each other. But to prove this is fairly complicated, so the proof is omitted.

*Proof of Proposition 4.1:* Assume first that $p \in (1/d, 1/(4\sqrt{d+1}))$. Since $G$ contains $\mathbb{T}_d$ as a subgraph it is clear that $p_c(G) \le 1/d$, so percolation occurs at $p$. Suppose for contradiction that $X$ contains a unique infinite cluster. Then by Lemma 4.2 the probability that two vertices are connected via paths of open edges is bounded away from 0. We will show that this is not the case for the present value of $p$.

Fix a vertex $t$ of $\mathbb{T}_d$. Then $(t,z)$ is a vertex of $G$ for every $z \in \mathbb{Z}$. We will prove that $P_p((t,0) \leftrightarrow (t,k))$ tends to 0 as $k$ increases, by means of the following path counting argument. In order for $\{(t,0) \leftrightarrow (t,k)\}$ to occur, there must be a path of open edges of some length $n \ge k$ between the two vertices. The number of paths from $(t,0)$ to $(t,k)$ of length $n$ is bounded by

$$4^n (d+1)^{n/2} ,$$

where the factor $4^n$ is a bound for the number of ways of specifying for each edge in the path whether it goes up or down in the $\mathbb{Z}$-direction or away from or towards $t$ in the $\mathbb{T}_d$-direction, and the $(d+1)^{n/2}$ is for specifying where the at most $n/2$ edges leading away from $t$ in $\mathbb{T}_d$ are heading. The probability that a given path of length $n$ is open is $p^n$, and so

$$P_p((t,0) \leftrightarrow (t,k)) \ \leq \ \sum_{n=k}^{\infty} (4p\sqrt{d+1})^n \,,$$

which tends to 0 as $k \to \infty$, as desired.

Next we move to the second part of the theorem, so assume now that $p > 1/2$. We will use Lemma 4.3. Let $u = (t_u, z_u)$ and $v = (t_v, z_v)$ be two arbitrary vertices of $G$. We need to show that for any $a > 0$ we can pick $N$ so large that $B(u,N)$ is connected to $B(v,N)$ with probability at least $1-a$. Pick a bi-infinite path $Q$ in $\mathbb{T}_d$ passing through $t_u$ and $t_v$. Then $H := Q \times \mathbb{Z}$ is a subgraph of $G$ isomorphic to $\mathbb{Z}^2$ and both $u$ and $v$ belong to $H$. Since $p > p_c(\mathbb{Z}^2)$ we know that the percolation restricted to $H$ will a.s. contain an infinite cluster and by Theorem 2.4 it will be unique. Put $A_v(N)$ and $A_v(N)$ for the events that $B(u,N)$ and $B(v,N)$ respectively intersect the infinite cluster of the percolation restricted to $H$. When $N$ is large enough $P_p(A_u(N)) = P_p(A_v(N)) > 1 - a/2$ and so by Bonferroni's inequality

$$P_p(B(u,N) \leftrightarrow B(v,N)) \geq P_p(A_u(N) \cap A_v(N)) > 1 - a \,,$$

and the proof is complete. □

An alternative proof of the $p > \frac{1}{2}$ part of Proposition 4.1, not using Lemma 4.3 but instead the notion of *cluster frequency* of Lyons and Schramm [82], will be outlined in Section 8. This part of Proposition 4.1 was generalized by Häggström et al. [57] to a result that states that if $G$ is the product of $d$ infinite connected graphs then $p_u(G) \leq p_c(\mathbb{Z}^d)$. The proof in [57] is a similar application of Lemma 4.3 as the one given here, but can, if one so wishes, be replaced by the cluster frequency technique in Section 8.

## 5. Uniqueness monotonicity

This section is devoted to the result that percolation on transitive graphs with $p > p_u$ a.s. produces a unique infinite cluster. First we need to define a few concepts and we begin with *Cayley graphs*.

Definition 5.1 *Let $H$ be a countable group and suppose that $S$ is a finite symmetric set of generators. Define a graph $G = G(H,S)$ with vertex set $H$ such that $\{u,v\}$ is an edge iff there exists an element $s \in S$ such that $u = vs$. The graph $G$ is called a (right)* **Cayley graph** *of $H$.*

Recall from Definition 1.2 the notion of automorphisms of a graph $G$. The set of all such automorphisms is denoted $\mathrm{Aut}(G)$, and constitutes a group under the operation of composition. A subgroup $H$ of $\mathrm{Aut}(G)$ is said to *act transitively* on $V$ if for all two vertices $u$ and $v$ there is an element $h \in H$ taking $u$ to $v$. All Cayley graphs $G$ are transitive since the underlying group $H$ itself can via left multiplication be identified with a subgroup of $\mathrm{Aut}(G)$ acting transitively on $G$. Most transitive graphs one comes across in the literature, including all the examples seen so far in this paper, are Cayley graphs. However not all transitive graphs are Cayley graphs; we will shortly give an example of such a graph.

A slightly larger class of transitive graphs is the class of *unimodular* transitive graphs, defined as follows.

Definition 5.2 *Let $G = (V, E)$ be a transitive graph and let $\mathrm{Aut}(G)$ denote the group of automorphisms of $G$. Let $H$ be a closed subgroup of $\mathrm{Aut}(G)$. For $v \in V$, define the* **stabilizer** *of $v$ with respect to $H$ as the subgroup $\mathrm{Stab}_H(v) = \{h \in H : hv = v\}$. Then $H$ is said to be* **unimodular** *if for every $u, v \in V$*

$$|\mathrm{Stab}_H(u)v| = |\mathrm{Stab}_H(v)u|.$$

*The graph $G$ is said to be* **unimodular** *if there exists a unimodular closed subgroup $H$ of $\mathrm{Aut}(G)$ that acts transitively on $V$.*

A Cayley graph $G$ of a countable group $H$ is always unimodular: it is trivial that $H$ itself represents a unimodular closed subgroup since $|\mathrm{Stab}_H(u)v| = 1$ for all $u$ and $v$. In fact the whole automorphism group $\mathrm{Aut}(G)$ is then also unimodular, see [10, Sect. 6]. On the other hand, there are unimodular transitive graphs that are not Cayley graphs: it is fairly easy to verify that all planar biregular[3] graphs are unimodular, while Chaboud and Kenyon [22] showed that a planar biregular graph with degree $d$ and dual degree $k$ is a Cayley graph iff $i|k$ for some $i \in \{2, 3, \ldots, d\}$. For instance, no planar biregular graph where $k$ is prime and $d < k$ is a Cayley graph.

The following simple example of a graph which is transitive but not unimodular (and therefore not a Cayley graph) is due to Trofimov [102]. In the binary tree $\mathbb{T}_2$, pick a directed bi-infinite path. For a vertex on the path, call its neighbor in the direction of the path its *parent*. From this one can in a natural way identify any pair of neighboring vertices of $\mathbb{T}_2$ as the parent and the child. Now create a new graph $G$ by adding, for every vertex, an edge between that vertex and its grandparent. Then $G$ is clearly transitive. However $\mathrm{Aut}(G)$ is not unimodular: If $u$ is a vertex and $v$ its parent, then $|\mathrm{Stab}(u)v| = 1$ whereas $|\mathrm{Stab}(v)u| = 2$.

Next we discuss the important *mass-transport principle*, which is due to Benjamini et al. [10]. Consider an automorphism invariant random process $X$ assigning values from some finite set $F$ to the edges of the graph $G = (V, E)$; automorphism invariance means that for any edge set $\{\{x_1, y_1\}, \ldots, \{x_k, y_k\}\}$ and any automorphism $f$ of $G$, the distribution of

$$(X(\{x_1, y_1\}), \ldots, X(\{x_k, y_k\}))$$

is the same as that of

$$(X(\{f(x_1), f(y_1)\}), \ldots, X(\{f(x_k), f(y_k)\}))\,. \tag{5.1}$$

[3] A planar graph is said to be *biregular* if it is regular and moreover all faces have the same number of edges.

(The basic example of such a process is when $X$ is an iid bond percolation.) Let $m : V \times V \times F^E \to \mathbb{R}_+$ be a function which is invariant under the diagonal action of $\mathrm{Aut}(G)$, i.e. such that for every automorphism $g$ and every $u, v \in V$ and every $\omega \in F^V$ one has $m(u, v, \omega) = m(gu, gv, g\omega)$. Intuitively we think of $m(u, v, \omega)$ as the amount of *mass* that is transported from $u$ to $v$ when $X = \omega$.

THEOREM 5.3 (THE MASS-TRANSPORT PRINCIPLE) *If $G$ is transitive and unimodular, then the expected amount of mass transported out of a vertex equals the expected amount of mass transported into it, i.e. for any vertex $u$*

$$\mathbb{E} \sum_{v \in V} m(u, v, X) = \mathbb{E} \sum_{v \in V} m(v, u, X).$$

In fact, the statement of the theorem holds in greater generality; it suffices that $X$ is invariant under the action of a closed transitive subgroup $H \subseteq \mathrm{Aut}(G)$. The mass-transport principle as formulated here fails on non-unimodular graphs, although a variant involving a certain "re-weighting" of mass does hold; see [10]. The proof of Theorem 5.3 gets considerably simpler if one restricts to the case of Cayley graphs, so we will settle for that; for a proof in the general case we again refer to [10].

*Proof of Theorem 5.3 for Cayley graphs:* This is straightforward: If $u$ and $v$ are vertices, then they are also members of the underlying group $H$ and there is a unique element $h = vu^{-1} \in H$ such that $v = hu$. Thus

$$\begin{aligned} \sum_{v \in V} \mathbb{E} m(u, v, X) &= \sum_{h \in H} \mathbb{E} m(u, hu, X) = \sum_{h \in H} \mathbb{E} m(h^{-1}u, u, h^{-1}X) \\ &= \sum_{h \in H} \mathbb{E} m(h^{-1}u, u, X) = \sum_{v \in V} \mathbb{E} m(v, u, X) \end{aligned}$$

where the third equality follows from the automorphism invariance of $X$. □

THEOREM 5.4 (HÄGGSTRÖM AND PERES [56]) *Let $G$ be a unimodular transitive graph. Consider iid bond percolation on $G$ and let $N$ denote the number of infinite clusters. If $p_1 < p_2$ and $P_{p_1}(N = 1) = 1$, then also $P_{p_2}(N = 1) = 1$.*

*Proof.* Let $p_2 > p_1 > p_c(G)$. Pick $X_{p_1}, X_{p_2} \in \{0, 1\}^E$ according to the simultaneous coupling construction (Coupling 1.1), so that in particular $X_{p_1}$ has distribution $P_{p_1}$ and $X_{p_2}$ has distribution $P_{p_2}$. The theorem follows if we can show that a.s. every infinite cluster of $X_{p_2}$ contains an infinite cluster of $X_{p_1}$. Since $X_{p_2}(e) \geq X_{p_1}(e)$ for every $e$, it suffices to show that every infinite $X_{p_2}$-cluster *intersects* an infinite $X_{p_1}$-cluster.

For each vertex $u$ define

$$D_1(u) = \inf\{\mathrm{dist}(u, v) : v \text{ is in an infinite cluster of } X_{p_1}\},$$

where $\mathrm{dist}(u, v)$ is the graph-theoretical distance between $u$ and $v$. Put $C(u, X_{p_i})$ for the connected component of $u$ in the open subgraph given by $X_{p_i}$. Let $A(u)$ be the event $\{D_1(u) = \min_{v \in C(u, X_{p_2})} D_1(v) > 0\}$. In words $A(u)$ is the event that $u$ is not in an infinite $X_{p_1}$-cluster, but no vertex in $u$'s $X_{p_2}$-component is closer to an infinite $X_{p_1}$-cluster than $u$ itself. Now if $X_{p_2}$ contains an infinite cluster that intersects no infinite cluster of $X_{p_1}$, then that cluster must contain a vertex $u$ closest to an infinite $X_{p_1}$-cluster and so for this vertex the event $B(u) := A(u) \cap \{|C(u, X_{p_2})| = \infty\}$ occurs. Thus we will be done as soon as we have established that

$$P(B(u)) = 0. \tag{5.2}$$

Write $M(u)$ for the number of vertices of $u$'s infinite $X_{p_2}$-cluster that are closest to an infinite $X_{p_1}$-cluster, and partition $B(u)$ as $B(u) = B^\infty(u) \cup B^f(u)$, where $B^\infty(u)$ is $B(u)$ intersected with the event that $M(u) = \infty$ and $B^f(u)$ is $B(u)$ intersected with the event that $M(u) < \infty$. Equation (5.2) will be verified when we have shown that $P(B^\infty(u)) = P(B^f(u)) = 0$. For $B^f(u)$ we use the mass-transport principle. On this event, let $m(x, y, X) = 1/M(x)$ if $y$ is a vertex in $C(x, X_{p_2})$ that is closest to an infinite cluster of $X_{p_1}$, and $m(x, y, X) = 0$ otherwise. Since the outgoing mass from any vertex is bounded by 1, the expected incoming mass to a vertex is also bounded by 1 by the mass-transport principle. However on $B^f(u)$ the incoming mass to some vertices will be infinite, and so we must have $P(B^f(u)) = 0$.

It remains to prove that $P(B^\infty(u)) = 0$. Write $B^\infty(u) = \cup_{k=1}^\infty B_k^\infty(u)$ where $B_k^\infty(u) = B^\infty(u) \cap \{D_1(u) = k\}$. That $B_k^\infty(u)$ has probability 0 can be seen by looking at the whole situation the following way: First condition on $X_{p_1}$ and then condition on $X_{p_2}(e)$ for all edges $e$ not incident to vertices within distance $k-1$ from infinite $X_{p_1}$-clusters. The conditional distribution on these remaining edges is then iid, with probability $(p_2 - p_1)/(1 - p_1)$ for an edge to be open in $X_{p_2}$. But if $B_k^\infty(u)$ occurs then infinitely many disjoint paths of length $k$ of such edges can tie the infinite $X_{p_2}$-cluster of $u$ to an infinite $X_{p_1}$-cluster and so this will a.s. happen. This proves that $P(B_k^\infty(u)) = 0$ for every $k$ and thereby also that $P(B^\infty(u)) = 0$. We have thus established (5.2), and the proof is complete. □

Our next task is to prove Schonmann's generalization of Theorem 5.4; see Theorem 5.6 below. The following lemma (Lemma 1.1 of [95]) states that an infinite cluster a.s. contains arbitrarily large balls.

Lemma 5.5 (Schonmann [95]) *Suppose that the iid bond percolation $X$ on a transitive graph $G = (V, E)$ a.s. contains an infinite cluster and let $M$ be an arbitrary positive integer. Then a.s. every infinite cluster of $X$ contains a ball of radius $M$.*

*Proof.* Fix a vertex $u$ and let $A = A(u)$ be the event that $u$ is in an infinite cluster that does not contain a ball of radius $M$. We want to prove that $P(A) = 0$. For all positive integers $k$, put $A_k$ for the event that $u$ is connected to a vertex at distance $2kM$ from $u$ and that the restriction of $X$ to the ball of radius $2kM$ centered at $u$ does not contain a ball of radius $M$ in its open subgraph. Clearly $A_k \downarrow A$ so that $P(A_k) \downarrow P(A)$. If it can be shown that for some constant $a = a(M) > 0$ we have $P(A_{k+1}|A_k) < 1 - a$, then it will thus follow that $P(A_k) \downarrow 0$ as desired.

Suppose that $\eta$ is a configuration in $A_{k+1}$. Then there is a vertex $x$ at distance $2kM + M$ from $u$ that is connected to $u$. Pick the first such $x$ (according to some arbitrary ranking of the vertices) and map $\eta$ to the configuration $\eta'$ which agrees with $\eta$ outside the ball of radius $M$ centered at $x$ and for which all edges of that ball are open. Clearly $\eta' \in A_k \setminus A_{k+1}$ and since $G$ has bounded degree there is a bounded number of configurations that are mapped to the same $\eta'$. Therefore $P(A_{k+1}^c|A_k) > a$ for some $a > 0$, as desired. □

In the proof of the general result we will also need the concept of "growing the cluster of a vertex $u$" which is an iterative process designed to find the cluster of $u$ by looking at the edges one by one in the following way: Suppose that $u$ is a vertex of a transitive graph $G$ on which a bond percolation process takes place. Order the edges of $G$ is some arbitrary way, put $C_0(u) = \{u\}$ and $\partial C_0(u) = \emptyset$. Then given $C_{n-1}(u)$ and $\partial C_{n-1}(u)$, check if there is an edge that goes between a vertex in $C_{n-1}(u)$ and $V \setminus C_{n-1}(u)$ that is not in $\partial C_{n-1}(u)$. If there is no such edge the process stops and we get $C(u) = C_{n-1}(u)$ and $\partial C(u) = \partial C_{n-1}(u)$. If there is such an edge, look at the smallest one, $e$, and check if it is open or closed in $X$. If $e$ is open then let $C_n(u)$ be $C_{n-1}(u)$ with the end vertex of $e$ outside $C_{n-1}(u)$ added and $\partial C_n(u) = \partial C_{n-1}(u) \cap \partial_E C_n(u)$. If $e$ is closed then $C_n(u) = C_{n-1}(u)$ and $\partial C_n(u) = \partial C_{n-1}(u) \cup \{e\}$. Note that if the cluster of $u$ is infinite then the process never terminates, whereas if the cluster of $u$ is finite the final $C_n(u)$ consists of the vertices in this cluster and $\partial C_n(u)$ coincides with its edge boundary.

THEOREM 5.6 (SCHONMANN [95]) *Let $G$ be a transitive graph. Consider iid bond percolation on $G$ and let $N$ denote the number of infinite clusters. If $p_1 < p_2$ and $P_{p_1}(N = 1) = 1$, then also $P_{p_2}(N = 1) = 1$.*

*Proof.* Use the same notation as in the proof of Theorem 5.4. The present situation differs from the one of Theorem 5.4 in that we do not assume that $G$ is unimodular and that we therefore do not have access to the mass-transport principle. However, since mass-transport was not used to show that $P(B^\infty(u)) = 0$, it suffices to prove that $P(B^f(u)) = 0$.

Put $B^f(u) = E(u) \cup E'(u)$ where $E'(u) = B^f(u) \cap \{D_1(u) \geq 2\}$ and $E(u) = B^f(u) \cap \{D_1(u) = 1\}$. We claim that if $P(E(u)) > 0$, then also $P(E'(u)) > 0$ and so it suffices to show that $P(E'(u)) = 0$. To prove the claim, assume that $P(E(u)) > 0$ and let $E^M(u)$ be the event that $E(u)$ occurs and that all the vertices of $u$'s infinite $X_{p_2}$-cluster that are incident to an infinite $X_{p_1}$-cluster are within distance $M$ from $u$. Since $E^M(u) \uparrow E(u)$, we have for $M$ large enough that $P(E^M(u)) > 0$. Now changing any $X$-configuration by changing the state of $(X_{p_1}(e), X_{p_2}(e))$ to $(0, 0)$ for every $e$ within distance $M$ from $u$ only changes its probability by a bounded constant factor; cf. Coupling 2.5. But this change maps any configuration in $E^M(u)$ to a configuration where $E'(v)$ happens for some vertex $v$ (in one of the infinite clusters that remain from $u$'s infinite $X_{p_2}$-cluster).

To prove that $P(E'(u)) = 0$ we fix an arbitrary $a > 0$ and show that $P(E'(u)) < a$. To this end, pick an integer $M$ so large that a given ball of radius $M$ intersects an infinite cluster of $X_{p_1}$ with probability exceeding $1-a$.

In order to make things work out smoothly we use what Schonmann refers to as a "duplication trick": Let $Z' = (Z'_{p_1}, Z'_{p_2})$ and $Z'' = (Z''_{p_1}, Z''_{p_2})$ be two independent $(\{0,1\}^E)^2$-valued random objects with the same distribution as $X = (X_{p_1}, X_{p_2})$. We start by growing the cluster of $u$ in $Z'_{p_2}$ as defined before the statement of the theorem, with the small adjustment that we stop the process if we find ourselves in a situation where $C_n(u)$ contains a ball of radius $M$. By Lemma 5.5 this will a.s. happen if the cluster of $u$ is infinite.

Put $I$ for the random time when the process stops. In words, $I$ is the first time when we have either found that the $Z'_{p_2}$-cluster of $u$ is finite and contains no ball of radius $M$ or when we have found that it contains such a ball. If we find ourselves in the latter situation we say that the event $F_1$ has occurred. It was just noted that $E'(u)$ is up to a set of measure 0 contained in $F_1$. On $F_1$, put $Y$ for the center of this ball. (If there is more than one possible such center, then pick one according to some arbitrary ranking of the vertices.) Define $Z = (Z_{p_1}, Z_{p_2}) \in (\{0,1\}^E)^2$ by putting $Z = Z'$ on all internal edges of $C_I(u)$ and all edges of $\partial C_I(u)$ and putting $Z = Z''$ on all other edges. Clearly $Z$ has the same distribution as $X$ so it will be sufficient to prove that $E'(u)$ cannot happen for $Z$. (Note that the way $Z$ is obtained from $Z''$ is similar to the local modifier (Coupling 2.5), except that here $Q$ is deterministic.)

Let $F_2$ be the sub-event of $F_1$ where the ball at $Y$ intersects an infinite $Z''_{p_1}$-cluster and note that on $F_2$ we have, by the construction of $Z$, that the $Z_{p_2}$-cluster of $u$ is within distance 1 from an infinite $Z_{p_1}$-cluster, in particular $E'(u)$ does not occur. We have shown

$$P(E'(u)) \leq P(F_1 \cap F_2^c).$$

Since $F_1$ is measurable with respect to $Z'$, and since $Z''$ is independent of $Z'$, we have by the choice of $M$ that

$$P(F_2|F_1) \geq 1-a.$$

It follows that $P(E'(u)) \leq P(F_1 \cap F_2^c) \leq P(F_2^c|F_1) < a$, as desired. □

Before moving on in the next two sections to the problem of determining when $p_c < p_u$, we end the present section with briefly mentioning a tantalizing open problem concerning the infinite clusters in the non-uniqueness regime $(p_c, p_u)$: that of so-called *cluster repulsion*. Suppose for some $p$ in the non-uniqueness regime of a transitive nonamenable graph $G$ that we can find two infinite clusters $C_1$ and $C_2$ that come within distance 1 from each other in infinitely many places. Call two such infinite clusters *strongly neighboring*. By the same reasoning as in the last paragraph of the proof of Theorem 5.4, it is easy to see that a.s. in the simultaneous coupling construction, $C_1$ and $C_2$ will have merged at level $p+\epsilon$ for any $\epsilon > 0$. It seems intuitively plausible that, already at level $p$, $C_1$ and $C_2$ ought to have merged, so that in other words the existence of two strongly neighboring infinite clusters has probability 0. If this is the case for any $p$, then $G$ is said to exhibit cluster repulsion. Häggström et al. [57] conjectured that cluster repulsion holds for any transitive graph, and supplemented their conjecture with a non-transitive counterexample. Very recently, Timár [101] proved the conjecture under the additional assumption of unimodularity, by means of an ingenious mass-transport argument. The nonunimodular case remains open.

## 6. The non-uniqueness phase for Cayley graphs

In this section and the next, we prove the two most general results known today on the existence of a non-uniqueness phase for nonamenable transitive graphs: The case of Cayley graphs with a suitably chosen set of generators of Pak and Smirnova-Nagnibeda [88] which is treated in this section, and the case of planar graphs with one end of Benjamini and Schramm [13], which we defer to Section 7. We begin with some preliminaries.

Let $G = G(H, S)$ be the nonamenable Cayley graph of the countable group $H$ with the symmetric set $S$ of generators. Let $d = |S|$ be the degree of $G$. Put $\phi = \phi(G) = \phi(H, S) = \kappa_E(G)/d$. Clearly, $\phi < 1$, and since $G$ is nonamenable $\phi > 0$. For two vertices $u$ and $v$, let $p^{(n)}(u, v)$ denote the probability that simple random walk on $G$ started at $u$ is at $v$ after $n$ steps. Put $\rho = \rho(G) = \rho(H, S) = \limsup_{n\to\infty}[p^{(n)}(u, v)]^{1/n}$, the **spectral radius** of $G$. This quantity is independent of $u$ and $v$, and since $G$ is nonamenable, $\rho < 1$, see [74]. In fact, the statement $\rho < 1$ is equivalent to nonamenability of $G$, and one can quantify the relation between $\phi$ and $\rho$, see [85]:

$$\phi \geq \frac{d(1-\rho)}{d-1}. \tag{6.1}$$

From (6.1) it follows that if $\{G_n\}$ is a sequence of Cayley graphs such that $\rho(G_n) \to 0$, then $\phi(G_n) \to 1$. This fact is all that we will need from (6.1).

The following two lemmas can be found in [12]:

LEMMA 6.1 *Let $G$ be a nonamenable graph with isoperimetric constant $\kappa_E$. Then*

$$p_c(G) \leq \frac{1}{\kappa_E + 1}.$$

*Proof.* Consider percolation with $p > 1/(1 + \kappa_E)$. Fix a vertex $u$ and grow the cluster of $u$ as described in the previous section. If $p < p_c$ then the process will a.s. eventually terminate with a finite $C_n(u)$ such that all edges of its edge boundary $\partial_E C_n(u)$ have been found to be closed. The number of open edges found in the process will then be $|C_n(u)|$ and the number of closed edges will be at least $|\partial_E C_n(u)|$. In other words, the fraction of the edges that will be found to be open will be at most

$$\frac{|C_n(u)|}{|\partial_E C_n(u)| + |C_n(u)|} \leq \frac{1}{1+\kappa_E}.$$

However, the states of the edges looked at in the process form an iid sequence $\{X_i\}$ where $P(X_i = 1) = p > 1/(1+\kappa_E)$ and so by standard large deviation theory there is a positive probability that there is no positive integer $n$ such that the fraction of the first $n$ $X_i$'s exceeds $1/(1+\kappa_E)$. In other words, there is a positive probability that the process of growing the cluster at $u$ never terminates, which is exactly what we wanted to prove. □

Lemma 6.2 *If $G$ is a transitive graph with degree $d$ for which*

$$\rho(G)p_c(G)d < 1,$$

*then $p_c(G) < p_u(G)$.*

*Proof.* The idea is the same as in the proof of Proposition 4.1: Let $p < 1/(d\rho)$ and show that $P_p(u \leftrightarrow v)$ becomes arbitrarily small when $u$ and $v$ are far enough apart. The result then follows from Lemma 4.2. Here we will again use a duplication trick: Consider two independent copies $X_1$ and $X_2$ of the percolation. If for some $\epsilon > 0$ we can find vertices $u$ and $v$ arbitrarily far apart such that $u$ is connected to $v$ in $X_1$ with probability at least $\epsilon$, then the event

$$B(u,v) := \{u \text{ is connected to } v \text{ in both } X_1 \text{ and } X_2\}$$

has probability at least $\epsilon^2$. Thus it suffices to prove that $P_p(B(u,v)) \to 0$ as $\text{dist}(u,v) \to \infty$.

Fix $u$ and let $a$ be a constant larger than but close enough to 1 so that $ap < 1/(d\rho)$. For $n \in \mathbb{Z}_+$ put $N_n(u)$ for the number of paths of length $n$ from $u$ to itself, i.e. $N_n(u) = p^{(n)}(u,u)d^n$. By definition of $\rho$, there exists an $n_0$ such that $n \geq n_0$ implies that $p^{(n)}(u,u) \leq (a\rho)^n$, whence $N_n(u) \leq (ad\rho)^n$. Fix a vertex $v$ at distance $k \geq n_0$ from $u$.

Now, for $B(u,v)$ to occur, it is necessary that there exists a self-avoiding path of edges open in $X_1$ from $u$ to $v$ and one such path of edges open in $X_2$. If we concatenate these two paths we get a path of some length $n \geq 2k$ from $u$ to itself that passes $v$ once, such that the edges before $v$ are open in $X_1$ and the edges after $v$ are open in $X_2$. For a given path from $u$ to itself passing $v$ once of length $n$, the probability that its edges are open in this way is $p^n$. Thus the expected number of such open paths is bounded by $N_n(u)p^n \leq (apd\rho)^n$. Consequently

$$P_p(B(u,v)) \leq \sum_{n=2k}^{\infty} (apd\rho)^n \to 0$$

as $k \to \infty$. □

We are now ready to deal with the Pak–Smirnova-Nagnibeda result. What we will show is that with any given symmetric set $S$ of generators, $p_c(H,S^k) < p_u(H,S^k)$ for $k$ large enough. Here $S^k$ is the multiset of elements of $H$ of the type $s_1 s_2 \ldots s_k$, $s_1, s_2, \ldots, s_k \in S$, where each element that can be so produced is taken with multiplicity equal to the number of ways that it can be written this way. Thus, to be entirely correct, $G_k := G(H,S^k)$ is a Cayley multigraph rather than a Cayley graph. Note that the degree of $G_k$ is $|S|^k = d^k$.

Theorem 6.3 (Pak and Smirnova-Nagnibeda [88]) *With $G_k$ as above and $k$ large enough,*

$$p_c(G_k) < p_u(G_k).$$

*Proof.* By the definition of $\rho(H, S) = \rho(G_1)$ the probability that simple random walk on $G_1$ started at $u$ is at $v$ after $n$ steps, decays like $e^{-\rho(G_1)n}$. Since $n$ steps of random walk on $G_k$ is the same as $kn$ steps of random walk on $G_1$ it thus follows that $\rho(G_k) \leq \rho(G_1)^k$; in particular $\rho(G_k) \to 0$. As observed above, this also entails that $\phi(G_k) \to 1$.

If we can show that $\rho(G_k)p_c(G_k)d^k \to 0$, then the result follows from Lemma 6.2. However, by Lemma 6.1,

$$\rho(G_k)p_c(G_k)d^k \leq \rho(G_k)\frac{1}{\kappa_E(G_k)}d^k = \frac{\rho(G_k)}{\phi(G_k)} \to 0.$$

□

## 7. The non-uniqueness phase for planar graphs

In this section we prove the result of Benjamini and Schramm [13] on the existence of a non-uniqueness regime for percolation on nonamenable transitive planar graphs:

Theorem 7.1 (Benjamini and Schramm [13]) *Let $G$ be a nonamenable planar transitive graph with one end. Then $p_c(G) < p_u(G)$.*

An important ingredient in the proof is planar duality, so this is a good place to recall some basic facts about planar duals.

In general, the planar dual $G^\dagger$ of an infinite planar graph $G$ is a multigraph and may have finite vertex set, but when $G$ is infinite with degree at least 3, as will always be the case here, $G^\dagger$ is an infinite graph.

Transitivity of $G$ does not guarantee the same for $G^\dagger$; in fact $G^\dagger$ may not even be regular as many of the standard planar lattices reveal. However, if $G^\dagger$ is regular (in which case one says that $G$ is *biregular*) then it is also transitive. This is a well-known fact, but the only written proof we are aware of can be found in [54]. More generally: $G$ is quasi-transitive if and only if $G^\dagger$ is quasi-transitive. (Recall that a graph is said to be quasi-transitive if its automorphism group partitions its vertex set into finitely many orbits, instead of only one as in the transitive case.)

The dual of a planar transitive graph $G$ may not have bounded degree (consider e.g. $G = \mathbb{T}_2$) but if $G$ is assumed to have one end, then $G^\dagger$ has bounded degree and also one end. This is so because if $G^\dagger$ has vertices of infinite degree, then every vertex $u$ of $G$ is incident to a face with infinitely many boundary edges, and exactly two neighbors of $u$ are incident to the same infinite face; clearly, the removal of $u$ plus a neighbor of $u$ that is not incident to the same infinite face partitions $G$ into at least two infinite subgraphs, contradicting that $G$ has one end. Thus since $G^\dagger$ is quasi-transitive it has bounded degree. Also, since $G^\dagger$ has bounded degree, any finite set $W$ of its vertices is connected to the rest of the graph by only finitely many edges. Thus $W$ is surrounded by a cycle of edges in $G$. This cycle is in turn surrounded by a cycle in $G^\dagger$. This means that the removal of $W$ does not partition $G^\dagger$ into more than one infinite subgraph and so $G^\dagger$ has one end.

An important observation is the following. Any finite connected set of vertices in $G$ is surrounded by a cycle in $G^\dagger$ and if $G$ has one end, then also every finite connected set of vertices in $G^\dagger$ is surrounded by a cycle in $G$.

The Benjamini–Schramm result that is the topic of this section is concerned with nonamenable transitive planar graphs with one end. By [13, Prop. 2.1] such graphs are always unimodular and, in fact, the conclusion holds with transitivity weakened to quasi-transitivity as shown in Lyons and Peres [80, Sect. 7]. We are thus free to use mass-transport ideas to prove the result. Our first step towards Theorem 7.1 is a generalization of Lemma 6.1 to percolation processes that are not necessarily iid but only automorphism invariant. Suppose that $X$ is an invariant bond percolation on a unimodular graph $G = (V, E)$ and assume that $X$ a.s. contains no infinite clusters. Put $p$ for the probability that a given edge is open. For all vertices $u$ let, as above, $C(u, X)$ be the cluster of $u$ in $X$. Define a mass-transport function $m$ by letting $m(u, v, X) = 1/(|C(v, X)| - |\partial_V C(v, X)|)$ if $u \in \partial_V C(v, X)$ and $v \notin \partial_V C(v, X)$. Otherwise let $m(u, v, X) = 0$. In other words all vertices that are in the boundary of a cluster distribute unit mass among all the "inner" vertices of that cluster. The probability that a given vertex is an inner vertex of a cluster is bounded from below by the probability the all edges incident to it are open, which in turn is, by Bonferroni's inequality, at least $1 - d(1 - p)$. Since $X$ contains only finite clusters a.s., the expected incoming mass to a vertex is thus at least

$$\frac{(1 - d(1 - p))\kappa_V}{1 - \kappa_V}.$$

On the other hand, the expected outgoing mass from a vertex is bounded from above, by the probability that the vertex is in the boundary of a cluster, which is in turn bounded by the probability that at least one edge incident to it is closed, i.e. by $d(1 - p)$. By the mass-transport principle we get

$$d(1 - p) \leq \frac{(1 - d(1 - p))\kappa_V}{1 - \kappa_V}$$

which entails

$$p \leq 1 - \frac{\kappa_V}{d}.$$

Consequently,

Lemma 7.2 (Benjamini et al. [10]) *Let $X$ be an automorphism invariant bond percolation on a unimodular transitive graph $G$ with degree $d$ such that the probability that a given edge is open exceeds $1 - \kappa_V(G)/d$. Then, with positive probability, $X$ contains infinite clusters.*

If $G$ is nonunimodular and $\mathrm{Aut}(G)$ is nonamenable, then Lemma 7.2 also holds, in the sense that there exists an $\eta > 0$ such that if edges are open with probability at least $1 - \eta$, then $X$ contains infinite clusters with positive probability. However, it is important to be aware that nonamenability of $G$ does *not* imply nonamenability of $\mathrm{Aut}(G)$. The converse, however, is true. Furthermore, for a unimodular graph, amenability of the graph itself and its automrophism group are equivalent. See [10] for more details, including the definition of amenability for groups.

We will need the following result, which is of course highly interesting also in its own right.

Theorem 7.3 (Benjamini et al. [9]) *Let $G = (V, E)$ be a unimodular nonamenable transitive graph and let $X$ be a critical iid bond percolation on $G$. Then $X$ a.s. contains no infinite clusters.*

Removing the unimodularity condition in this result seems at present like a tough nut to crack, although see Peres et al. [91] and Timár [100] for some progress in this direction.

*Proof of Theorem 7.3.* First we rule out the possibility of having a unique infinite cluster. Assume that $X$ a.s. contains a unique infinite cluster $U$. For $v \in V$ let $U(v)$ be the set of vertices $u$ of $U$ such that $\mathrm{dist}(v, u) = \mathrm{dist}(v, U)$, i.e. the set of vertices of $U$ that are closest to $v$.

Let $a > 0$ be a small number and let $Y_a$ be an iid bond percolation on $G$ with retention parameter $a$, independent of $X$. Put $X \setminus Y_a$ for the percolation process one gets by taking an edge to be open if it is open in $X$ and closed in $Y_a$. Then $X \setminus Y_a$ is iid percolation with parameter $p_c(G) - a$ and thus sub-critical.

Now define another percolation process $Z_a$ by declaring each edge $e = \{u, v\} \in V$ to be open in $Z_a$ if $\mathrm{dist}(u, U) < 1/a$, $\mathrm{dist}(v, U) < 1/a$ and $U(u) \cup U(v)$ is contained in a connected component of $X \setminus Y_a$. Then a.s.

$$\lim_{a \downarrow 0} P(Z_a(e) = 1 | X) = 1\,,$$

so by the Dominated Convergence Theorem

$$\lim_{a \downarrow 0} P(Z_a(e) = 1) = 1.$$

Thus $a$ may be picked so small that $P(Z_a(e) = 1) > 1 - \kappa_V(G)/d$, so by Lemma 7.2, $Z_a$ with positive probability contains an infinite cluster. However, if $v_1, v_2, v_3, \ldots$ are the vertices of an infinite self-avoiding open path in $Z$ and $u_i \in U(v_i)$ then by the definition of $Z_a$ all the $u_i$'s are in the same cluster of $X \setminus Y_a$, and the sequence $\{u_i\}$ contains infinitely many distinct vertices. This contradicts the sub-criticality of $X \setminus Y_a$.

Now assume that $X$ contains infinitely many infinite clusters. Recall from the proof of Theorem 2.4 the definition of a trifurcation: a trifurcation is a vertex such that exactly three of the edges incident to it are open and the removal of it splits an infinite cluster into three infinite clusters. In the proof of Theorem 2.4 it was argued that the set $T$ of trifurcations is with positive probability nonempty iff $X$ contains infinitely many infinite clusters a.s. It follows from ergodicity (cf. the proof of Lemma 2.6) that indeed $T$ will a.s. be nonempty in this case. Furthermore we claim that if a trifurcation $v$ is removed from the open subgraph of $X$, then each of the three infinite clusters that this results in, contains another trifurcation. To prove this define a mass-transport function $m$ by putting $m(u, v, X) = 1$ if $v$ is the unique trifurcation that is closest in the open subgraph of $X$ to $u$, and $m(u, v, X) = 0$ otherwise. Then the expected outgoing mass from a vertex $v$ is bounded by 1. However if $v$ is a trifurcation where one of the three infinite clusters its removal results in has no other trifurcation, then $v$ would get unit mass from all vertices in that cluster. Thus, a.s., this cannot happen.

We shall now construct a forest (i.e. a graph whose connected components are trees) $F = F(X)$ with $T$ as vertex set in such a way that the distribution of $F$ is automorphism-invariant: First assign to each $v \in T$ a random variable $\xi(v)$ that is uniformly distributed on $[0, 1]$ and independent of $X$ and $\{\xi(u) : u \neq v\}$. Then for a given $v \in T$, let $X_1$, $X_2$ and $X_3$ denote the three infinite sub-clusters of $X$ that the removal of $v$ results in. For each $i$ put an edge between $v$ and the one among the trifurcations, $t$, in $X_i$ that are closest in $X$ to $v$ that has the smallest $\xi(t)$. The graph $F$ constructed in this way is clearly automorphism-invariant and it is also a forest. To see the latter, assume for contradiction that $F$ contains a cycle $(v_0, v_1, v_2, \ldots, v_k)$, $v_0 = v_k$. Then by the construction of $F$ we must have that $\text{dist}(v_i, v_{i+1})$ is the same for all $i$. Therefore $\xi(v_0) > \xi(v_1) > \ldots > \xi(v_k) = \xi(v_0)$, a contradiction.

Note that every component of $F$ is contained in an infinite cluster of $X$.

Again we will use an iid bond percolation $Y_a$ as above, independent of $X$. Define a percolation $Z_a$ on $F$ by letting an edge be open if its two end-vertices are in the same cluster of $X \setminus Y_a$. Since $X \setminus Y_a$ contains no infinite clusters, neither does $Z_a$. For $v \in T$ let $K(v)$ denote the cluster of $v$ in $Z_a$ and let $\partial_F K(v)$ denote $K(v)$'s vertex-boundary in $F$. Since the components of $F$ are binary trees, at least half of the vertices of $K(v)$ are in $\partial_F K(v)$. Therefore, for a given vertex $v \in V$,

$$P(v \in T, v \in \partial_F K(v)) \geq \frac{1}{2} P(v \in T).$$

However, the probability that an edge in $F$ is closed in $Z_a$ tends to 0 as $a$ tends to 0, whence the left-hand side tends to 0. On the other hand, the right-hand side is nonzero and independent of $a$, a contradiction. □

Note that for the first part of the proof (no *unique* infinite cluster at criticality), the fact that $X$ is iid percolation was only used through its property of automorphism invariance. Thus, the same proof yields the following more general result, which will serve as one more ingredient in the proof of Theorem 7.1.

COROLLARY 7.4 *Let $X$ be an automorphism invariant bond percolation on a unimodular transitive nonamenable graph such that $X$ a.s. has a unique infinite cluster, $C$. Then a.s. $p_c(C) < 1$.*

The following lemma is from [13, Lem. 3.3].

LEMMA 7.5 *Let $X$ be an iid bond percolation on a nonamenable planar transitive graph $G$ with one end and let $X^\dagger$ be the dual percolation on $G^\dagger$. Then a.s. at least one of $X$ and $X^\dagger$ contains infinite clusters.*

*Proof.* If $X$ as well as $X^\dagger$ contain only finite clusters a.s. then, since $G$ and $G^\dagger$ both have one end, every cluster of $X$ is surrounded by a circuit in $X^\dagger$ and vice versa. Define the **rank** of a cluster of $X$ as follows. First say that a cluster has rank 0 if it does not contain a circuit that surrounds a cluster of $X^\dagger$ that in turn contains a circuit that surrounds a cluster of $X$. Then recursively set the rank of a cluster $C$ of $X$ as one plus the maximum rank of clusters that are surrounded by a circuit of a cluster of $X^\dagger$ that is surrounded by a circuit of $C$. Define the rank of a vertex $v$ to be the rank of its cluster.

Now define an automorphism-invariant percolation $Y$ on $G$ by saying that an edge is open if the rank of neither of its end-vertices exceeds a given number $R$. Clearly the probability that an edge is open tends to 1 as $R \to \infty$, so we can let $R$ be chosen so that this probability exceeds $1 - \kappa_V(G)/d$. Then by Lemma 7.2, $Y$ must contain infinite clusters. On the other hand, since every cluster in $X$ is surrounded by clusters of arbitrarily large rank, $Y$ cannot contain infinite clusters, and we have a contradiction. □

THEOREM 7.6 *Let $G$ be a nonamenable planar transitive graph with one end, let $X$ be iid bond percolation on $G$ and let $X^\dagger$ be the dual percolation on $G^\dagger$. Then a.s. either $X$ and $X^\dagger$ both contain infinitely many infinite clusters or one of them contains a unique infinite cluster and the other contains no infinite clusters.*

*Proof.* By Lemma 7.5 at least one of $X$ and $X^\dagger$ must contain infinite clusters. If $X$ contains infinitely many infinite clusters, these must be separated by some infinite cluster of $X^\dagger$ and vice versa. Also if $X$ a.s. contains infinitely many infinite clusters, then pick a finite connected subgraph $H$ of $G$ so large that $H$ with positive probability intersects two infinite clusters in $G \setminus H$. Since there is also a positive probability that all edges of $H$ are open in $X$ we find that with positive probability $X^\dagger$ contains at least two infinite clusters. Hence if $X$ contains infinitely many infinite clusters, the so does $X^\dagger$ and vice versa.

It remains to rule out the possibility that $X$ and $X^\dagger$ both contain a unique infinite cluster. In order to obtain a contradiction, consider the graph $H$ formed by letting the vertices of $H$ be the union of $V$, $V^\dagger$ and the points where an edge of $G$ crosses an edge of $V^\dagger$ (considering a fixed given embedding of $G$ in the plane). The edges of $H$ are the "half egdes" of $G$ and $G^\dagger$, i.e. formally the pairs $\{c, v\}$ where $c$ is a crossing and $v$ is one of the end vertices of one of the two edges that cross at $c$. Note that $H$ is quasi-transitive. Define an invariant percolation, $Y$, on $H$ by declaring an edge $\{c, v\}$ to be open iff $\{u, v\}$ is open (in $X$ or $X^\dagger$), where $u$ is the other end-vertex (in $G$ or $G^\dagger$) of the edge in $G$ or $G^\dagger$ that has a crossing in $c$. Now if $X$ and $X^\dagger$ both have a unique infinite cluster a.s. then $Y$ has a.s. exactly two infinite clusters. It is now tempting to refer to Lemma 2.6 to make the conclusion that this is impossible. However, $Y$ is not an iid percolation and therefore we cannot do so. Luckily we have other means to deal with this situation: Assume that $Y$ a.s has two infinite clusters. Pick one of the clusters, $C$, uniformly at random, and put $W$ for the set of edges of $H$ connecting $V(C)$ to the vertices of the component of $G \setminus C$ containing the other infinite cluster. Put $W^\dagger$ for the set of edges in the dual of $H$ crossing edges of $W$. Then $W^\dagger$ is an invariant percolation on the dual of $H$ and by planarity it consists of a infinite path. Hence $p_c(W^\dagger) = 1$, a contradiction to Corollary 7.4. □

*Proof of Theorem 7.1:* First we claim that $p_u(G) = 1 - p_c(G^\dagger)$: Let $X$ be percolation on $G$ with retention parameter $p$ and and let $X^\dagger$ be the dual percolation. If $p < p_u(G)$ then by Theorem 7.6, $X^\dagger$ contains (infinitely many) infinite clusters and so $1 - p \geq p_c(G^\dagger)$. On the other hand if $p > p_u(G)$ then the same theorem tells us that $X^\dagger$ does not contain infinite clusters and so $1-p \leq p_c(G^\dagger)$.

Next we claim that $p_c(G) + p_c(G^\dagger) < 1$ for if this had not been the case then with $p = p_c(G)$, by Theorem 7.3, none of $X$ or $X^\dagger$ contains an infinite cluster, another contradiction to Theorem 7.6.

Putting this together yields

$$p_c(G) < 1 - p_c(G^\dagger) = p_u(G)$$

as desired. □

## 8. Uniqueness for $p$ close to 1

When is $p_u < 1$? On the tree $\mathbb{T}_d$, $d \geq 2$, one has $p_u = 1$ even though $p_c < 1$. This phenomenon also occurs on any transitive graph $G$ with more than one end for which $p_c(G) < 1$, since the removal of a finite set of edges splits the graph into at least two infinite connected graphs and there is always a positive probability that all edges of that set are closed. So in order to have any hopes for proving uniqueness for some $p < 1$, one-endedness needs to be assumed. We shall start by returning to the $p > \frac{1}{2}$ half of Proposition 4.1 in the Grimmett–Newman example $\mathbb{T}_d \times \mathbb{Z}$ (which is obviously one-ended). Here we restate that half (with a weaker condition on $d$, but note that the $d \geq 17$ condition there was not necessary for the $p \geq \frac{1}{2}$ part):

Proposition 8.1 *Let $d \geq 2$, and let $G = \mathbb{T}_d \times \mathbb{Z}$. For iid bond percolation on $G$ with $p > \frac{1}{2}$, we get a.s. a unique infinite cluster.*

Let us sketch an alternative proof of this result based on the notion of *cluster frequency*, invented by Lyons and Schramm [82], rather than on Lemma 4.3. (An ancestor of cluster frequency for the $\mathbb{Z}^d$ setting is *cluster density*, for which an analog of Theorem 8.4 below was established by Burton and Keane [19].) Let $G = (V, E)$ be a transitive graph, and let $Z(0), Z(1), \ldots$ be a simple random walk on $G$ starting from some fixed vertex. If now $Y$ is a $\{0, 1\}^V$-valued random object satisfying automorphism invariance (the obvious analog of (5.1) with bonds replaced by sites), then $Y(Z(0)), Y(Z(1)), \ldots$ becomes a stationary process. Thus, the ergodic theorem tells us that the limit

$$\lim_{n\to\infty} \frac{1}{n} \sum_{i=0}^{n-1} Y(Z(i)) \tag{8.1}$$

exists. Furthermore, we have the following.

LEMMA 8.2 *Let $G = (V, E)$, $Z(0), Z(1), \ldots$ and $Y$ be as above. Then the limit in* (8.1) *depends a.s. only on $Y$ (and not on the random walk).*

*Proof.* Fix a $c \in [0, 1]$, and write $L$ for the limit in (8.1). By Lévy's 0-1-law, we have a.s. that

$$\lim_{m\to\infty} P(L \le c \,|\, Z(0), \ldots, Z(m), Y) \;=\; I_{\{L\le c\}}\,.$$

Hence, for any $\epsilon > 0$, we can find an $M < \infty$ such that with probability at least $1 - \epsilon$ we have

$$P(L \le c \,|\, Z(0), \ldots, Z(M), Y) \;\in\; [0, \epsilon] \cup [1 - \epsilon, 1]\,. \tag{8.2}$$

We clearly have

$$\lim_{n\to\infty} \frac{1}{n} \sum_{i=M}^{M+n-1} Y(Z(i)) \;=\; L\,,$$

and substituting in (8.2) gives

$$P\left(\lim_{n\to\infty} \frac{1}{n} \sum_{i=M}^{M+n-1} Y(Z(i)) \le c \;\middle|\; Z(0), \ldots, Z(M), Y\right) \;\in\; [0, \epsilon] \cup [1 - \epsilon, 1]\,,$$

where the left-hand side depends on $Z(0), \ldots, Z(M)$ only via $Z(M)$, so that

$$P\left(\lim_{n\to\infty} \frac{1}{n} \sum_{i=M}^{M+n-1} Y(Z(i)) \le c \;\middle|\; Z(M), Y\right) \;\in\; [0, \epsilon] \cup [1 - \epsilon, 1]\,. \tag{8.3}$$

Due to automorphism invariance of $Y$ in conjuncion with $Z(M)$ being independent of $Y$, we have that the left-hand side in (8.3) has the same distribution as $P(L \le c \,|\, Y)$. Since $\epsilon > 0$ was arbitrary, this means that $P(L \le c \,|\, Y) \in \{0, 1\}$ a.s., and since $c \in [0, 1]$ was arbitrary the proof is complete. □

Lemma 8.2 motivates the following definition.

DEFINITION 8.3 *Let $G = (V, E)$ be a transitive graph, and fix $y \in V$ as well as a bond percolation configuration $\xi \in \{0, 1\}^E$, and write $C_y$ for the cluster of*

*$\xi$ containing $y$. Let $Z(0), Z(1), \ldots$ be simple random walk on $G$ starting from some fixed vertex. The* **cluster frequency** *of $C_y$ is defined as*

$$\lim_{n\to\infty} \frac{1}{n}\sum_{i=0}^{n-1} I_{\{Z(i)\in C_y\}}\,,$$

*provided that the limit exists a.s. and is independent of the random walk.*

Lyons and Schramm [82] established the following.

THEOREM 8.4 *For iid bond percolation*[4] *on a transitive graph $G = (V, E)$, we get with probability 1 that every cluster has a well-defined cluster frequency.*

*Proof.* We shall employ a coupling idea that yields a different, and in our opinion even more instructive, proof than the one in [82]. Fix $y \in V$, pick the bond percolation process $X \in \{0,1\}^E$ according to $P_p$, and pick $Y \in \{0,1\}^V$ as follows: for each cluster $C$ of $X$ independently, let all vertices of $C$ take value 0, or let them take value 1, with probability 1/2 each. The distribution of $Y$ is obviously automorphism invariant, and we get that

$$\lim_{n\to\infty} \frac{1}{n}\sum_{i=0}^{n-1} Y(Z(i)) \qquad (8.4)$$

exists, and by Lemma 8.2 it is independent of the random walk.

Next, pick $Y' \in \{0,1\}^V$ by setting

$$Y'(x) = \begin{cases} 1 - Y(x) & \text{if } x \in C_y \\ Y(x) & \text{otherwise,} \end{cases}$$

and note that $Y'$ has exactly the same distribution as $Y$. Hence, also

$$\lim_{n\to\infty} \frac{1}{n}\sum_{i=0}^{n-1} Y'(Z(i)) \qquad (8.5)$$

exists and is independent of the random walk. Combining what we know about the limits in (8.4) and (8.5), we get that

$$\lim_{n\to\infty} \frac{1}{n}\sum_{i=0}^{n-1} (Y(Z(i)) - Y'(Z(i))) \qquad (8.6)$$

exists and is independent of the random walk, and the same thing for

$$\lim_{n\to\infty} \frac{1}{n}\sum_{i=0}^{n-1} |Y(Z(i)) - Y'(Z(i))| \qquad (8.7)$$

[4] We focus here on the iid setting, but in fact the result as well as its proof go through without change in the more general setting of automorphism invariant percolation.

because all terms in (8.6) have the same sign (positive if $Y(y) = 1$ and negative otherwise). But

$$|Y(Z(i)) - Y'(Z(i))| \;=\; I_{\{Z(i) \in C_y\}}$$

so that the limit in (8.7) is in fact the cluster frequency of $C_y$. $\Box$

*Proof of Proposition 8.1:* Write $\theta_{\mathbb{Z}^2}(p)$ for the probability in iid bond percolation on $\mathbb{Z}^2$ that the origin is in an infinite cluster, and recall from Section 1 that $\theta_{\mathbb{Z}^2}(p) > 0$ when $p > \frac{1}{2}$. By the same token as in the proof of Lemma 4.2, we get that

$$P_p(u \leftrightarrow v) \;\geq\; \theta_{\mathbb{Z}^2}(p)^2 \tag{8.8}$$

for any $u, v \in \mathbb{Z}^2$.

Switching to percolation on $\mathbb{T}_d \times \mathbb{Z}$, recall from Section 4 that for any two vertices $u$ and $v$ in this graph there is a subgraph of $\mathbb{T}_d \times \mathbb{Z}$ isomorphic to $\mathbb{Z}^2$, and that we can therefore conclude that (8.8) holds for any $u, v \in \mathbb{T}_d \times \mathbb{Z}$ as well.

Fix a vertex $y$ in $\mathbb{T}_d \times \mathbb{Z}$ and then pick a $p' \in (\frac{1}{2}, p)$. It follows from (8.8) that with positive $P_{p'}$-probability, the cluster $C_y$ has a cluster frequency $\mathrm{freq}(C_y)$ satisfying

$$\begin{aligned} \mathrm{freq}(C_y) \;&\geq\; \theta_{\mathbb{Z}^2}(p')^2 \\ &>\; 0\,, \end{aligned}$$

and note that if so happens, then $C_y$ must be infinite. But the existence of *some* infinite cluster $C$ with $\mathrm{freq}(C) > 0$ is a tail event, so that by Kolmogorov's 0-1-law

$$P_{p'}(\exists \text{ an infinite cluster } C \text{ with } \mathrm{freq}(C) > 0) \;=\; 1\,. \tag{8.9}$$

The sum of $\mathrm{freq}(C)$ over all clusters $C$ is easily seen to be bounded by 1. Hence, using (8.9), we get with $P_{p'}$-probability 1 that there exist finitely many infinite clusters maximizing cluster frequency (among all clusters in the particular percolation configuration that we happened to get); call such a frequency-maximizing cluster *special*.

Now consider the usual enhancement from $p'$ to $p$, so that in other words $X_{p'}, X_p \in \{0,1\}^E$ are picked jointly as in Coupling 1.1. We saw in Theorem 5.4 that a.s. any infinite cluster in $X_p$ will contain some infinite cluster in $X_{p'}$. But an inspection of the proof of Theorem 5.4 shows that it goes through virtually unchanged to show that any infinite cluster in $X_p$ contains a *special* infinite cluster in $X_{p'}$. Hence $X_p$ has a.s. only finitely many infinite clusters, and by Lemma 2.6 it therefore has only one. $\Box$

We next move on from the Grimmett–Newman example to more general one-ended transitive graphs and a result of Babson and Benjamini [6] in the direction of establishing a uniqueness phase for $p$ close to 1. First we note that when $G$ is planar with one end it was shown in the proof of Theorem 7.1 that $p_u(G) = 1 - p_c(G^\dagger)$ and since $G^\dagger$ has bounded degree, $p_c(G^\dagger) > 0$ and hence $p_u(G) < 1$. Thus the planar case is already settled. But in fact the assumption of planarity can be removed, provided that the group is taken to be finitely presented:

Theorem 8.5 (Babson and Benjamini [6]) *Let $G = (V, E)$ be a one-ended nonamenable Cayley graph of a finitely presented group. Then $p_u(G) < 1$.*

The assumption of nonamenability is not essential in any other respect than that it ensures that $p_c(G) < 1$. The result also trivially holds for amenable graphs with one end provided that $p_c(G) < 1$.

We will be content with sketching some of the ideas behind Theorem 8.5. Define a *cutset* to be a set $F$ of edges such that with $F$ removed the graph $(V, E \setminus F)$ has at least two infinite connected components. Note that since $G$ has one end a cutset must be infinite. If $X$ is an iid bond percolation on $G$ with more than one infinite cluster, then the set of closed edges $X$ must contain a cutset, so if we can show that for $p$ large enough the probability for this is 0, then we will be done.

For that we use the fact that there exists a number $D < \infty$ such that a cutset on $G$ must contain an infinite set $\{f_1, f_2, f_3, \ldots\}$ of distinct edges such that $\text{dist}(f_i, f_{i+1}) \leq D$ for all $i$, where the distance between two edges is the minimal distance between two of their end-vertices. This fact is a direct consequence of [6, Cor. 4] whose proof we omit; a simpler proof by Timár (also omitted here) can be found in [99] and in [80, Sect. 6.6].

Now let $d$ be the degree of $G$. The probability that there exists set of edges $\{f_1, f_2, \ldots\}$ of the kind just described with the first $n$ of the $f_i$'s closed in $X$, is bounded by $(2d^D(1-p))^n$. Hence the probability that $X$ contains a cutset of closed edges is 0 when $p \geq 1 - 1/(2d^D)$, implying Theorem 8.5.

## 9. The situation at $p_u$

The question in focus in this section is what happens when $p = p_u$. In the case of percolation on planar graphs with one end the answer is more or less already given by the results of Section 7:

Theorem 9.1 (Benjamini and Schramm [13]) *Let $G = (V, E)$ be a planar unimodular transitive graph with one end. Then percolation at $p_u(G)$ a.s. produces a unique infinite cluster.*

*Proof.* Let $X$ be an iid bond percolation on $G$ with $p = p_u(G)$ and let as usual $X^\dagger$ denote the dual percolation on $G^\dagger$. By the proof of Theorem 7.1, $p_u(G) + p_c(G^\dagger) = 1$ so $X^\dagger$ is critical percolation on $G^\dagger$. The graph $G^\dagger$ is unimodular (this was discussed in Section 7) and we may therefore apply Theorem 7.3 to deduce that $X^\dagger$ contains no infinite clusters. Hence, by Theorem 7.6, $X$ contains a unique infinite cluster. □

The behavior on planar graphs is not general. The next result of Schonmann shows that for the Grimmett–Newman example of Section 4 the situation is the opposite.

Theorem 9.2 (Schonmann [94]) *Let $G = \mathbb{T}_d \times \mathbb{Z}$. Then a.s. percolation at $p_u(G)$ does not produce a unique infinite cluster. In particular if $d$ is chosen so that $p_c < p_u$, e.g. if $d$ is sufficiently large, then percolation at $p_u$ a.s. produces infinitely many infinite clusters.*

Schonmann's result was soon afterwards generalized by Peres to the following, while Lyons and Schramm [82] proved that uniqueness of the infinite cluster at $p_u$ is a.s. impossible for Caley graphs of so-called Kazhdan groups.

Theorem 9.3 (Peres [90]) *Let $G = (V_G, E_G)$ and $H = (V_H, E_H)$ be two infinite connected transitive graphs and assume that $G$ is nonamenable and unimodular. Then iid bond percolation on $G \times H$ at $p_u(G \times H)$ a.s. does not produce a unique infinite cluster.*

Actually, as consequence of the remark after Lemma 7.2, the condition that $G$ is nonamenable and unimodular can be replaced by the weaker assumption that $\mathrm{Aut}(G)$ is nonamenable. Recall, however, that when $G$ is nonunimodular, this does not automically follow from $G$ itself being nonamenable.

Since Theorem 9.3 is more general and has a proof that is not harder to follow than the original proof of Theorem 9.2, we present only the proof of Theorem 9.3. Note that Theorem 9.3 implies that $p_u < 1$, a special case of the result from [57] mentioned at the end of Section 5.

*Proof of Theorem 9.3.* Assume that $p_0$ is such that percolation on $G \times H$ with retention parameter $p_0$ a.s. produces a unique infinite cluster. We will show that then there exists $p_* < p_0$ such that percolation with retention $p_*$ also produces a unique infinite cluster. We will do this by showing that for some $p_* < p_0$ it is the case that for any $a > 0$ one can pick $N$ so large that $P_{p_*}(B(u,N) \leftrightarrow B(v,N)) > 1 - a$ for any $u, v \in V_G \times V_H$. The theorem then follows from Lemma 4.3.

We use the simultaneous coupling construction (Coupling 1.1) of percolation processes for all $p$: Assign to the edges, $e$, independent uniform [0,1] random variables $U(e)$ and denote the underlying probability measure by $P$. For $p \in [0,1]$ put $O(p)$ for the set of edges $e$ with $U(e) \leq p$ and note that $O(p)$ has the distribution of the open edges of iid percolation with retention $p$. Put $C(u,p)$ for the connected component of the vertex $u$ in $(V_G \times V_H, O(p))$ and for a set $W$ of vertices put $C(W,p) = \cup_{u \in W} C(u,p)$.

Since $G$ is nonamenable there exists, by Lemma 7.2, a number $\eta > 0$ such that any invariant site percolation process on $G$, for which a given vertex, is open with probability at least $1 - \eta$, with positive probability produces infinite clusters. Fix such an $\eta$.

Now put $C_\infty(p_0)$ for the unique infinite cluster in $(V_G \times V_H, O(p_0))$ and define

$$A_1 = A_1(r) = \{u \in V_G \times V_H : B(u,r) \cap C_\infty(p_0) \neq \emptyset\}.$$

Fix $r$ so large that $P(u \in A_1(r)) \geq 1 - \eta/6$. Next define

$$A_2 = A_2(r,n) = \{u \in V_G \times V_H : \forall v, w \in B(u, r+1) \cap C_\infty(p_0) : \mathrm{dist}(v, w; O(p_0)) < n\}.$$

Here $\text{dist}(v, w; O(p))$ refers to the distance between $v$ and $w$ in $(V_G \times V_H, O(p))$. Now fix $n$ so large that $P(u \in A_2(r, n)) \geq 1 - \eta/6$.

Let $d = d_G + d_H$ denote the degree of $G \times H$ and fix $p_1$, $p_2$ and $p_*$ so that

$$p_0 - \frac{\eta}{6d^{r+n}} < p_1 < p_2 < p_* < p_0.$$

Now define

$$A_3 = A_3(r, n, p_1) = \{u \in V_G \times V_H : U(e) \notin [p_1, p_0] \text{ for all edges } e \text{ in } B(u, r+n)\}.$$

Then $P(u \in A_3) \geq 1 - \eta/6$. Put $A = A_1 \cap A_2 \cap A_3$ and note that $P(u \in A) \geq 1 - \eta/2$. Fix two vertices $h_1, h_2 \in V_H$ and let

$$B = \{g \in V_G : (g, h_1) \in A \wedge (g, h_2) \in A\}.$$

Then $B$ is the set of open vertices of an invariant site percolation on $G$ such that $P(g \in B) \geq 1 - \eta$ for any given vertex $g$ and hence $B$, with positive probability, has infinite clusters. Since the event that $B$ contains infinite clusters is automorphism invariant and determined by the iid random variables $U(e)$, $B$ in fact contains an infinite clusters with probability 1.

Consider an infinite path $(g_1, g_2, g_3, \ldots)$ in $B$. Since $(g_j, h_1) \in A_1$, there is a vertex $u_j \in C_\infty(p_0)$ in $B((g_j, h_1), r)$. Since $(g_j, h_1) \in A_2 \cap A_3$ for all $j$, there is a path from $u_j$ to $u_{j+1}$ of length at most $n$ in $O(p_1)$. Concatenating these paths for $j = 1, 2, 3, \ldots$ gives an infinite path in $O(p_1)$ that comes within distance $r$ from $(g_j, h_1)$ for every $j$. By an analogous argument there is an infinite path $(v_1, v_2, \ldots)$ in $O(p_1)$ that comes within distance $r$ from $(g_j, h_2)$ for every $j$. Thus for any vertex $g \in V_G$ we have

$$P(B((g, h_1), r) \leftrightarrow B((g, h_2), r) \text{ in } O(p_2)) \mid |C(g, B)| = \infty) = 1. \tag{9.1}$$

This is so because there is an infinite number of edge-disjoint paths that can join the two paths above in $O(p_1)$ and all of these have a positive probability of being open at level $p_2$ whatever their status at level $p_1$. Similar arguments were used in Section 4.

Now pick $R_0$ so large that the probability that for $g \in V_G$ the ball $B(g, R_0)$ in $G$ intersects an infinite cluster of $B$ is at least $1 - a/2$, and put $R = r + R_0$. Taking (9.1) into account together with the triangle inequality shows that

$$P(B((g, h_1), R) \leftrightarrow B((g, h_2), R) \text{ in } O(p_2)) \geq 1 - \frac{a}{2}. \tag{9.2}$$

Fix two vertices $u = (g_u, h_u)$ and $v = (g_v, h_v)$ of $G \times H$. For $h \in V_H$ put $F_h$ for the event that $B(u, R) \leftrightarrow B((g_v, h), R)$ and $B(v, R) \leftrightarrow B((g_v, h, R)$. By (9.2), $F_h$ has probability at least $1 - a$ for any $h$. Therefore, $F_h$ occurs for infinitely many $h$ with probability at least $1 - a$. However, on this event the sets $C(B(u, R), p_2)$ and $C(B(v, R), p_2)$ come infinitely often within distance $2R + \text{dist}(g_u, g_v; E_G)$ of each other. Thus, since $p_* > p_2$,

$$P(B(u, R) \leftrightarrow B(v, R) \text{ in } C(p_*)) \geq 1 - a$$

as desired. □

## 10. Simultaneous uniqueness

Consider the simultaneous coupling construction (Coupling 1.1) of iid bond percolation on $\mathbb{Z}^d$ for all $p \in [0,1]$. We know from the results in Section 2 (Theorem 2.4, say) that for each $p > p_c$

$$P(X_p \text{ has a unique infinite cluster}) = 1 .$$

Since there are uncountably many such $p$'s, it is not a priori obvious whether or not the order of the quantifiers can be exchanged in this statement, i.e., whether

$$P(\text{for all } p > p_c,\ X_p \text{ has a unique infinite cluster}) = 1 . \tag{10.1}$$

This question of *simultaneous uniqueness* was first posed by Matthew Penrose, and appeared in print in Grimmett [35], where it was pointed out that (10.1) holds for $d = 2$. This follows easily from the contour argument in the proof of Theorem 2.2.

Soon thereafter, Alexander [4] established (10.1) for $d \geq 3$ as well as for a wider class of percolation processes in $d$-dimensional Euclidean space.

The issue of simultaneous uniqueness in the more general context of transitive graphs was later considered by Häggström and Peres [56] and Häggström et al. [57]. In that setting, (10.1) of course needs to be replaced by

$$P(\text{for all } p > p_u,\ X_p \text{ has a unique infinite cluster}) = 1 .$$

This was established under the unimodularity assumption in [56] and without that condition in [57]. In fact, the results in [56] and [57] have implications for the non-uniqueness phase $(p_c, p_u)$ as well. The general result is the following.

Theorem 10.1 *Consider the simultaneous coupling construction of iid bond percolation for all $p \in [0,1]$ on a transitive graph. With probability 1, we have for all $p_1, p_2$ such that $p_c < p_1 < p_2$, that any infinite cluster in $X_{p_2}$ contains an infinite cluster in $X_{p_1}$.*

In other words, when we gradually raise $p$ in the simultaneous coupling construction, the only thing that ever happens to infinite clusters after $p = p_c$ is that they grow and merge – no new infinite clusters are ever born.

To prove Theorem 10.1 is a matter of modifying the proof of Theorem 5.6 in such a way that the order in which edges are searched in the procedure exploring a cluster is based not on a pre-fixed ordering but on the uniform $[0,1]$ variables of Coupling 1.1, turning in effect the exploration procedure to so-called invasion percolation [23]. See [57] for the details.

## 11. Dynamical percolation

Similar in spirit to the question of simultaneous uniqueness considered in the previous section, is the study of so-called *dynamical percolation* (Häggström et al. [58]). Here the percolation process evolves in continuous time in the following manner. Fix $G = (V, E)$ and the parameter $p \in (0, 1)$. Each edge $e \in E$ independently flips between states 0 (closed) and 1 (open) according to a two-state continuous time Markov chain whose flip rate is $p$ when in state 0, and $1-p$ when in state 1. The unique stationary distribution for this Markov chain places mass $1-p$ on state 0 and mass $p$ on state 1. Hence, if we take the entire $\{0,1\}^E$-valued process $(X(t))_{t\geq 0}$ to begin with $X(0)$ chosen according to product measure $P_p$ on $\{0,1\}^E$, then $X(t)$ has distribution $P_p$ for all $t$. Thus, for any fixed $t \geq 0$,

$$P(X(t) \text{ has an infinite cluster}) = \begin{cases} 0 & \text{if } p < p_c \\ 1 & \text{if } p > p_c \,. \end{cases}$$

Analogously to the simultaneous uniqueness issue in Section 10, we have a nontrivial task if we wish to strengthen this statement by interchanging the order of the quantifiers, i.e., by saying that

$$P(\text{for all } t, X(t) \text{ has no infinite cluster}) = 1 \;\text{ if } p < p_c \tag{11.1}$$

and

$$P(\text{for all } t, X(t) \text{ has an infinite cluster}) = 1 \;\text{ if } p > p_c \,. \tag{11.2}$$

It turns out (see [58]) that (11.1) and (11.2) hold for dynamical percolation on arbitrary graphs, whereas the situation at $p = p_c$ is considerably more intricate: there exist graphs $G$ such that at the critical value there is $P_{p_c}$-a.s. no infinite cluster, while in the dynamical percolation process an infinite cluster shows up in $X(t)$ at some random exceptional times $t$; likewise there are graphs $G$ which $P_{p_c}$-a.s. have an infinite cluster, yet in the corresponding dynamical percolation process there is at certain exceptional times no infinite cluster. See Schramm and Steif [96] for a recent breakthrough amounting to that critical dynamical site percolation on the two-dimensional triangular lattice exhibits exceptional times.

Concerning uniqueness of the infinite cluster for supercritical percolation on $\mathbb{Z}^d$, Peres and Steif [92] showed that this can be strengthened to a "for all $t$" statement in dynamical percolation. More precisely, for $G = \mathbb{Z}^d$ and $p > p_c$, they demonstrated that the dynamical percolation process $(X(t))_{t\geq 0}$ satisfies

$$P(\text{for all } t, X(t) \text{ has a unique infinite cluster}) = 1 \,.$$

The same paper offers a detailed and beautiful analysis of the somewhat exotic ways in which the number of infinite clusters may vary in time in critical dynamical percolation on certain non-homogeneous trees.

## 12. Dependent percolation

Generally speaking, one of the most obvious (and, indeed, most-studied) directions for extending the iid bond percolation model studied in previous sections is to relax the iid assumption and allow for (various kinds of) dependence between edges. Such processes can be of great interest in themselves (as will be amply exemplified in the next section), or sometimes (as we saw in the proof of Theorem 7.3) even serve as tools in the analysis of iid percolation.

To work with *arbitrary* probability measures on $\{0,1\}^E$ would, however, be to take things a bit too far, as not much of interest can be said in such a general setting. So we still need *some* conditions.

For percolation on $\mathbb{Z}^d$, *translation invariance* is an extremely natural (and, in the percolation literature, almost universally employed) condition. A probability measure $\nu$ on $\{0,1\}^E$ (where $E$ is the edge set of the $\mathbb{Z}^d$ lattice) is said to be translation invariant if a $\{0,1\}^E$-valued random object with distribution $\nu$ has the following property: for any finite edge set $\{\{x_1,y_1\},\ldots,\{x_k,y_k\}\} \subset E$ and any $z \in \mathbb{Z}^d$, the distribution of

$$(X(\{x_1,y_1\}),\ldots,X(\{x_k,y_k\}))$$

is the same as that of

$$(X(\{x_1+z,y_1+z\}),\ldots,X(\{x_k+z,y_k+z\}))$$

(cf. the definition (5.1) of automorphism invariance).

Translation invariant percolation on $\mathbb{Z}^d$ does not in general inherit the uniqueness of the infinite cluster property of iid percolation. For a trivial example, consider the percolation process on $\mathbb{Z}^2$ where all vertical edges are a.s. open, while all horizontal edges are a.s. closed; clearly, this results in an infinite number of infinite clusters. A more challenging task is to construct, given a finite $k \geq 2$, a translation invariant percolation process resulting a.s. in precisely $k$ infinite clusters. But it can be done, as shown by Burton and Keane in the very interesting follow-up paper [20] to their famous breakthrough [19] discussed in Sections 1 and 2.

So more conditions are needed in order to deduce uniqueness of the infinite cluster. A striking sufficient condition (in the setting of translation invariance) is the notion of *finite energy*, dating back to Newman and Schulman [86]:

DEFINITION 12.1 *For an arbitrary graph $G=(V,E)$, a probability measure $\nu$ on $\{0,1\}^E$ is said to have* **finite energy**[5] *if it admits conditional probabilities such that for a $\{0,1\}^E$-valued random object $X$ with distribution $\nu$, any $e \in E$, and any $\xi \in \{0,1\}^{E\setminus\{e\}}$, we have*

$$0 < \nu(X(e)=1 \mid X(E\setminus\{e\})=\xi) < 1\,. \tag{12.1}$$

In fact, the uniqueness result of Burton and Keane [19] (earlier quoted as Theorem 2.4) was originally stated in the setting of finite energy percolation processes on $\mathbb{Z}^d$:

[5]Lyons and Schramm [82] found it fruitful to distinguish between *insertion tolerance* (meaning that the first inequality in (12.1) holds), and *deletion tolerance* (ditto for the second inequality). If the probability in (12.1) is bounded away from 0 and/or 1, then we speak of *uniform finite energy*, *uniform insertion tolerance* and *uniform deletion tolerance*.

Theorem 12.2 (Burton and Keane [19]) *Any finite energy translation invariant bond percolation process on $\mathbb{Z}^d$ produces a.s. at most one infinite cluster.*

To adapt the proof of Theorem 2.4 to this non-iid setting, let us first inspect where in the proof of Theorem 2.4 the iid assumption is actually used. This happens in two places, namely,

(a) to deduce in Lemma 2.6 that the number of infinite clusters is an a.s. constant, and
(b) in the application of the local modifier (Coupling 2.5) to produce trifurcations.

Now, the a.s. constance of the number of infinite clusters in (a) is no longer true under the conditions of Theorem 12.2, but this is easily taken care of by ergodic decomposition of the probability measure in question. (Here one should verify that the finite energy property is preserved under ergodic decomposition, a detail that Burton and Keane slipped in [19] but which has been clarified several times later, e.g., in [32, Lem. 1] and [82, Lem. 3.6].)

As to (b), consider the following generalization of Coupling 2.5 to the non-iid setting.

Coupling 12.3 The local modifier in a non-iid setting. *Let $G = (V, E)$ be an infinite graph, and let $\nu$ be a probability measure on $\{0, 1\}^E$. A coupling of two $\{0, 1\}^E$-valued random objects $X$ and $X'$, both with distribution $\nu$, can be obtained as follows. Let $Q$ be a probability distribution on the set $\mathcal{E}$ of all finite subsets of $E$ such that $Q$ assigns positive probability to each element of $\mathcal{E}$.*

1. *Pick $F \in \mathcal{E}$ according to $Q$.*
2. *Pick $X(E \setminus F)$ according to the distribution prescribed by $\nu$, and set $X'(E \setminus F) = X(E \setminus F)$.*
3. *Pick $X(F)$ according to the conditional distribution that $\nu$ prescribes, given $X(E \setminus F)$.*
4. *Independently of Step 3, pick $X'(F)$ according to the conditional distribution that $\nu$ prescribes, given $X'(E \setminus F)$.*

Note that if $\nu$ equals product measure $P_p$, then Coupling 12.3 reduces to Coupling 2.5. We may furthermore note in the more general setting where $\nu$ has finite energy, all outcomes of $X(F)$ and $X'(F)$ have positive conditional probability given Steps 1 and 2. Therefore, we may for the construction of trifurcations in the proof of Theorem 12.2 proceed just as in the proof of Theorem 2.4, only replacing Coupling 2.5 by Coupling 12.3. The rest of the proof is identical to that of Theorem 2.4, so Theorem 12.2 is established.

## 13. Dependent percolation: examples

As elsewhere in mathematics, extending percolation-theoretic results to more general settings would be of questionable value unless these settings are shown to include some interesting concrete examples. For the theory of dependent percolation that we began to discuss in the previous section, there is no shortage of such examples.

One of the most important examples is the *random-cluster model*, introduced in a 1972 paper by Fortuin and Kasteleyn [31]. After a relatively calm period, interest in the model caught on again around 1990 following contributions by Edwards and Sokal [28], which clarified the model's intimate relation to Ising and Potts models, and Aizenman et al. [1], which initiated a modern probabilistic analysis of the model using stochastic domination and correlation inequality techniques. The level of activity has since then remained high, with a vast number of papers on applications to Ising and Potts models as well as studies of the random-cluster model as an object worthy of attention in its own right. This work has been surveyed in, e.g., Georgii et al. [34], Grimmett [38], and the recent book by Grimmett [40].

On a finite graph $G = (V, E)$, the random-cluster measure $\phi^G_{p,q}$ with parameters $p \in [0, 1]$ and $q > 0$ is defined as the probability measure on $\{0, 1\}^E$ which to each $\xi \in \{0, 1\}^E$ assigns probability

$$\frac{q^{k(\xi)}}{Z^G_{p,q}} \prod_{e \in E} p^{\xi(e)}(1-p)^{1-\xi(e)}\,, \tag{13.1}$$

where $k(\xi)$ is the number of connected components (including isolated vertices) of the graph $G_\xi = (V, \{e \in E : \xi(e) = 1\})$, and $Z^G_{p,q}$ is a normalizing constant. Note that taking $q = 1$ reduces to iid percolation with retention parameter $p$, whereas taking $q \neq 1$ in general leads to dependence between edges. Single-edge conditional distributions are easy to compute from (13.1): if $X \in \{0, 1\}^E$ is taken to have distribution $\phi^G_{p,q}$, then we have for any edge $e = \{x, y\} \in E$ and any $\xi \in \{0, 1\}^{E\setminus\{e\}}$ that

$$\phi^G_{p,q}(X(e) = 1 \,|\, X(E \setminus \{e\}) = \xi) = \begin{cases} p & \text{if } x \leftrightarrow y \text{ in } \xi \\ \frac{p}{p+(1-p)q} & \text{otherwise,} \end{cases} \tag{13.2}$$

where "$x \leftrightarrow y$ in $\xi$" means that there exists a path from $x$ to $y$ only using edges $e$ with $\xi(e) = 1$.

It is not immediately clear how one should extend the notion of random-cluster measures to the case where $G = (V, E)$ is infinite; note that $\{0, 1\}^E$ is then uncountable, so that each $\xi \in \{0, 1\}^E$ ought to receive probability 0, and a direct definition like (13.1) clearly will not do. Instead, two approaches have been proposed, namely (a) using DLR (Dobrushin–Lanford–Ruelle) equations, and (b) taking weak limits of finite-graph random-cluster measures.

In approach (a), originally considered in the 1960's for Gibbs measures in infinite volume, the idea is that single-edge conditional probabilities (and more generally conditional distributions on finite sets) should behave as on finite graphs. Thus, for an infinite graph $G = (V, E)$, we call a probability measure $\phi$ on $\{0, 1\}^E$ a *DLR random-cluster measure for $G$ with parameters $p$ and $q$* if (13.2) holds for all $e \in E$ and $\phi$-almost all $\xi \in \{0, 1\}^{E\setminus\{e\}}$. The corresponding agreement of conditional distributions on arbitrary finite edge sets follows easily (see, e.g., [34, Lem. 6.18]).

In approach (b), we first define a sequence $(G_n)_{n\geq 1}$ of finite subgraphs of $G$ such that

(i) $V_1 \subseteq V_2 \subseteq \cdots$,
(ii) each $x \in V$ is in all but at most finitely many $V_n$'s,
(iii) each $e = \{x, y\} \in E$ is in $E_n$ if and only if $x$ and $y$ are both in $V_n$,

and call such a sequence an *exhaustion* of $G$. Consider now the sequence $(\phi_{p,q}^{G_n})_{n\geq 1}$ of random-cluster measures for the graphs $(G_n)_{n\geq 1}$. It may happen that these converge weakly (in the sense of convergence of probabilities of all cylinder events) to a probability measure $\phi$ on $\{0,1\}^E$, in which case we call $\phi$ a *limiting random-cluster measure for $G$ with parameters $p$ and $q$*. A standard compactness argument shows that such a limit can always be obtained by passing to a subsequence if necessary.

For $q \geq 1$ more can be said. Then the right-hand side of (13.2) becomes increasing in $\xi$, an observation that opens up the door to a wonderful machinery of stochastic domination (much of it based on the so-called Holley Inequality, see, e.g., [34] for a detailed introduction to how this works) which among other things shows that, for any $n$, the projections on $\{0,1\}^{E_n}$ of the measures $(\phi_{p,q}^{G_{n+m}})_{m\geq 0}$ is stochastically increasing in $m$. It follows that the sequence $\phi_{p,q}^{G_1}, \phi_{p,q}^{G_2}, \ldots$ converges (without passing to a subsequence) to a probability measure on $\{0,1\}^E$ which thereby qualifies as a limiting random-cluster measure for $G$ with parameters $p$ and $q$; we denote it by

$$\phi_{p,q}^{G,free}$$

where "free" is short for "with free boundary condition", a notion that will be clarified through contrast to its "wired" counterpart in the next paragraph. Another consequence of the stochastic domination machinery is that $\phi_{p,q}^{G,free}$ is independent of the choice of exhaustion $(G_n)_{n\geq 1}$. For the case where $G$ is the $\mathbb{Z}^d$ lattice, we therefore get that $\phi_{p,q}^{G,free}$ is translation invariant.

An alternative limiting procedure (still when $q \geq 1$) is to replace $\phi_{p,q}^{G_1}, \phi_{p,q}^{G_2}, \ldots$ by the measures $\phi_{p,q}^{G_1,wired}, \phi_{p,q}^{G_2,wired}, \ldots$ defined by modifying $\phi_{p,q}^{G_n}$ in such a way that $k(\xi)$ in (13.1) counts only those connected components that do not intersect the boundary $\partial_V V_n$ (recall that $\partial_V V_n$ is defined as the set of all $x \in V_n$ that have at least one neighbor in $V \setminus V_n$). This is tantamount to viewing all vertices of the boundary $\partial_V V_n$ as a single one, or in other words to "wiring" them to each other. The same stochastic domination arguments as those alluded to in the previous paragraph (except that the stochastic monotonicity goes in the other direction) yield convergence to a probability measure

$$\phi_{p,q}^{G,wired}$$

where "wired" is short for "with wired boundary condition".

So far we have been silent on the issue of existence of DLR random-cluster measures. Avoiding the issue for just a few more lines, let us remark that it is straightforward to show that (still for $q \geq 1$) any DLR random-cluster measure for $G$ with parameters $p$ and $q$ must be sandwiched between $\phi_{p,q}^{G,free}$ and $\phi_{p,q}^{G,wired}$ in the sense of stochastic domination.

For the case where $G$ equals the $\mathbb{Z}^d$ lattice (or more generally an amenable transitive graph), Grimmett [36], and independently Pfister and Vande Velde [93], found the following elegant solution to the problem of existence of DLR random-cluster measures. The first key observation is that the single-edge conditional probabilities in $\phi_{p,q}^{G_n}$ are bounded away from 0 and 1 uniformly in $n$ (by (13.2), they are in fact always between $\frac{p}{p+(1-p)q}$ and $p$). This can be shown to imply the finite energy property of the limiting measure $\phi_{p,q}^{G,free}$, to which Theorem 12.2 therefore can be applied to deduce the $\phi_{p,q}^{G,free}$-a.s. uniqueness of the infinite cluster. Considering for a moment what it means to have an infinite cluster, we see that for any $\{x, y\} \in E$ we can $\phi_{p,q}^{G,free}$-a.s. find some sufficiently large finite subgraph of $G$ containing $x$ and $y$, such that by looking in the box we can verify either that $x$ and $y$ are in the same cluster, or that (at least) one of them is in a finite cluster not equal to the cluster containing the other. This last property is a kind of a.s. continuity in $\xi$ of the event $x \leftrightarrow y$ showing up in (13.2), and this turns out to be sufficient for establishing the DLR equation (13.2) for $\phi_{p,q}^{G,free}$. The same argument applies to $\phi_{p,q}^{G,wired}$, so we have the following.

Theorem 13.1 (Grimmett [36], Pfister and Vande Velde [93]) *For the random-cluster model on an amenable transitive graph $G$ with parameters $p \in (0,1)$ and $q \geq 1$, we have that $\phi_{p,q}^{G,free}$ and $\phi_{p,q}^{G,wired}$ are both DLR random-cluster measures, and that they both put full measure on the event of having at most one infinite cluster.*

Combining this result with the aforementioned sandwiching relation, shows that uniqueness of DLR random-cluster measures for an amenable transitive graph $G$ with parameters $p$ and $q \geq 1$ is equivalent to having

$$\phi_{p,q}^{G,free} \,=\, \phi_{p,q}^{G,wired}\,.$$

This last equality is known to hold, loosely speaking, for "most" but not all pairs $(p, q)$. For such a result, the amenability condition is essential; see Häggström [46] and Jonasson [71]. See also Grimmett [36, 38, 40] for more on this fascinating problem in the $\mathbb{Z}^d$ setting.

So what about DLR properties of $\phi_{p,q}^{G,free}$ and $\phi_{p,q}^{G,wired}$ on nonamenable graphs, where the uniqueness device Theorem 12.2 is not available? The problem was settled to the following unexpected answer by Jonasson [71] and Georgii et al. [34]: For $q \geq 1$ and an arbitrary graph $G$, $\phi_{p,q}^{G,free}$ is always a DLR random-cluster measure, whereas $\phi_{p,q}^{G,wired}$ need not be. Instead, $\phi_{p,q}^{G,wired}$ satisfies the modified DLR equation

$$\phi_{p,q}^{G}(X(e) = 1 \,|\, X(E \setminus \{e\}) = \xi) = \begin{cases} p & \text{if } x \overset{C}{\longleftrightarrow} y \text{ in } \xi \\ \frac{p}{p+(1-p)q} & \text{otherwise,} \end{cases} \tag{13.3}$$

where "$x \overset{C}{\longleftrightarrow} y$ in $\xi$" is short for the event that either $x \leftrightarrow y$ in $\xi$, or $x$ and $y$ are both in infinite clusters.

We may think of the definition of $x \overset{C}{\longleftrightarrow} y$ as a kind of "compactification of $G$ at infinity" of the graph: to walk from one vertex $x$ to another vertex $y$, it is allowed to walk from $x$ "to infinity", and then back from infinity to $y$. This alternative definition of DLR random-cluster measures goes back to Häggström [46], who studied such measures on trees, where the older choice (13.2) of DLR equation reduces, somewhat boringly, to iid percolation. Note that for graphs exhibiting $\phi_{p,q}^{G,wired}$-a.s. uniqueness of the infinite cluster, the right-hand sides of (13.2) and (13.3) coincide, so that this later development does not conflict with Theorem 13.1. Examples of graphs where $\phi_{p,q}^{G,wired}$ fails to exhibit uniqueness of the infinite cluster are studied in [46] and in Häggström et al. [54].

As a further application of the aforementioned stochastic domination machinery, we have for fixed $q \geq 1$ and any infinite $G$ that $\phi_{p,q}^{G,free}$ is stochastically increasing in $p$, and likewise for $\phi_{p,q}^{G,wired}$. The measures can be shown to satisfy a 0-1-law for the existence of infinite clusters, whence there exist critical values $0 \leq p_{c,wired} \leq p_{c,free} \leq 1$ (depending on $G$ and $q$) such that

$$\phi_{p,q}^{G,free}(\exists \text{ an infinite cluster}) = \begin{cases} 0 & \text{if } p < p_{c,free} \\ 1 & \text{if } p > p_{c,free} \end{cases} \tag{13.4}$$

and

$$\phi_{p,q}^{G,wired}(\exists \text{ an infinite cluster}) = \begin{cases} 0 & \text{if } p < p_{c,wired} \\ 1 & \text{if } p > p_{c,wired}\,. \end{cases} \tag{13.5}$$

For transitive graphs, Lyons [79] proved a random-cluster analog of the uniqueness monotonicity results in Section 5, thereby establishing the existence of critical values $0 \leq p_{u,wired} \leq p_{u,free} \leq 1$ satisfying analogs of (13.4) and (13.5) for the existence of a *unique* infinite cluster. We thus have four critical values $p_{c,wired}$, $p_{c,free}$, $p_{u,wired}$ and $p_{u,free}$ (depending on $G$ and $q$) for the random-cluster model on a transitive graph. As to the ordering of these values, there is − besides determining when the inequalities are strict − only one nontrivial case, namely the ordering between $p_{c,free}$ and $p_{u,wired}$. These problems are treated in some detail in [54], where among other things it is shown that both $p_{c,free} < p_{u,wired}$ and $p_{u,wired} < p_{c,free}$ can happen.

Let us now return to the more familiar situation where $G$ is the $\mathbb{Z}^d$ lattice. Even in this case, however, much less is understood about the random-cluster model for $q < 1$ compared to $q \geq 1$, due to the absence of various stochastic monotonicities. One interesting case, however, which nevertheless is fairly well-understood, is when $p$ and $q$ are sent to 0 simultaneously in such a way that $\frac{q}{p} \to 0$. Then, as shown by Häggström [45], any sequence of random-cluster measures (either in sense (a) or sense (b) above) converges weakly to the so-called *uniform spanning tree measure* of Pemantle [89]. This measure arises also (and is in fact defined) by the following limiting procedure: Let $(G_1, G_2, \ldots)$ be an exhaustion[6] of $\mathbb{Z}^d$, let $\mu_n$ denote the measure on $\{0,1\}^{E_n}$ corresponding to picking a spanning tree for $G_n$ at random, uniformly. The limiting measure

[6] In addition to conditions (i)–(iii) above, we now also need to take each $G_n$ to be connected.

$\mu = \lim_{n\to\infty} \mu_n$ on $\{0,1\}^E$ exists in the usual sense of convergence of cylinder probabilities.

Inspecting the properties of $\mu$, we find that it is translation invariant, but that it fails to satisfy the finite energy condition. Hence, the usual argument for uniqueness of the infinite cluster no longer works. But this couldn't possibly matter, we would like the reader to interject at this point, because surely a limit of spanning trees for $G_n$ will be a spanning tree for $\mathbb{Z}^d$ (and therefore connected)? Well, not always. In the beautiful main result of [89], it is shown that $\mu$ yields a.s. a spanning tree when $d \leq 4$, whereas for $d \geq 5$ an infinite collection of infinite trees is a.s. obtained. See Benjamini et al. [11, 8] for some subsequent fascinating work on the topic of these uniform spanning tree measures, on $\mathbb{Z}^d$ as well as on more general graphs.

On the topic of limits of random spanning trees as percolation-theoretic objects, we would like to mention that the other canonical way to choose a spanning tree for a finite graph at random – so-called *minimal spanning trees* – allows the same kind of infinite-volume limit. For a finite graph $G = (V, E)$, the minimal spanning tree arises by first assigning iid weights to the edges according to some continuous distribution (the precise choice of this distribution is inconsequential, so one usually takes it to be uniform on $[0, 1]$ as in Coupling 1.1), and then picking the spanning tree for $G$ whose total edge-weight is minimized (note that there will be no tie a.s.). The limit along an exhaustion $(G_1, G_2, \ldots)$ of $\mathbb{Z}^d$ yields a measure on $\{0,1\}^E$ that is closely related to invasion percolation. It is natural to ask whether connectivity in this limiting object depends on the dimension $d$ in similar fashion as the corresponding uniform spanning tree objects. Newman and Stein [86] gave heuristic arguments for getting a single spanning tree for $\mathbb{Z}^d$ when $d \leq 7$, and an infinite collection of such trees when $d \geq 9$ (the case $d = 8$ appears to be "the critical dimension" and even harder to decide), but so far this has been rigorously established only for $d = 2$ (Alexander [5]). See also Lyons et al. [81] for a study of these objects beyond the $\mathbb{Z}^d$ setting.

## 14. Dependent percolation on $\mathbb{Z}^2$

In this section, we will continue our study of dependent percolation, but restrict to the $\mathbb{Z}^2$ lattice, which has some particular geometric features not exhibited in transitive graphs more generally. The key feature appears to be the combination of planarity and amenability, and the arguments of this section go through *mutatis mutandis* for other (reasonable) planar lattices in two-dimensional Euclidean space.

We switch in this section to *site percolation* rather than bond percolation, thus considering a $\{0,1\}^{\mathbb{Z}^2}$-valued random object $X$. Similarly to bond percolation, we are interested in clusters of 1's, or 1-*clusters*, for short, meaning connected components in the subgraph of $\mathbb{Z}^2$ obtained by removing all vertices $x$ such that $X(x) = 0$ and all edges incident to such vertices. The main novel feature of the present section (compared to previous sections) is that we will simultaneously be interested in 0-*clusters*, defined in the same way as 1-clusters but with the roles of states 0 and 1 interchanged.

Rather than ask, as usual, for uniqueness of the infinite 1-cluster, we can be even more ambitious and ask for uniqueness of the infinite cluster in the setting where both 1-clusters and 0-clusters are considered simultaneously.

The uniqueness results of Section 2 (in particular Theorem 2.4) are straightforward to adapt to the setting of site percolation. Thus, for iid site percolation on $\mathbb{Z}^2$, we know that there is a.s. at most one infinite 1-cluster, and likewise a.s. at most one infinite 0-cluster. The question at hand is therefore whether or not an infinite 1-cluster and an infinite 0-cluster can coexist.

Another known fact of iid site percolation on $\mathbb{Z}^2$ is that the critical value $p_c$ for emergence of an infinite connected component exceeds $\frac{1}{2}$ (the strict inequality was first obtained by Higuchi [63], while the weaker statement $p_c \geq \frac{1}{2}$ is ancient – viewed on the time-scale of the history of percolation – and follows by combining the results of Harris [62] and Hammersley [60]). It follows that for no value of the retention parameter $p$ do an infinite 1-cluster and an infinite 0-cluster coexist.

We may now note that such a strong uniqueness result does not extend to the $\mathbb{Z}^3$ setting; this follows from the result of Campanino and Russo [21] that the critical value for site percolation on $\mathbb{Z}^3$ is strictly less than $\frac{1}{2}$, so that at the symmetry point $p = \frac{1}{2}$ (as well as in a nontrivial interval around it) we have a.s. coexistence of an infinite 1-cluster and an infinite 0-cluster. Coexistence of infinite clusters of both categories can happen also for planar lattices in the hyperbolic plane, as can easily be deduced from the results in Section 7. (However, Corollary 7.4 can be used to show that in the hyperbolic plane, coexistence of a *unique* infinite 1-cluster and a *unique* infinite 0-cluster a.s. does not happen, analogously to Theorem 7.6.)

On the other hand, the $\mathbb{Z}^2$ lattice does feel a bit narrow for encompassing a unique infinite 1-cluster together with a unique infinite 0-cluster, and results on the impossibility of such coexistence on $\mathbb{Z}^2$ ought to be available beyond the iid setting. An early result in this direction was obtained by Cognilio et al. [26] in the setting of Gibbs measures for the ferromagnetic Ising model. For more general results, conditions will be needed, and translation invariance alone will not do as shown, e.g., by counterexamples of Burton and Keane [19]. In the seminal paper by Gandolfi et al. [33], the notion of positive correlations was identified as being of central importance in this setting. Recall the following definition.

Definition 14.1 *For a finite or countable set $V$, a function $f : \{0,1\}^V \to \mathbb{R}$ is said to be* **increasing** *if it is increasing in each of its coordinates. A probability measure $\nu$ on $\{0,1\}^V$ is said to exhibit* **positive correlations** *if for any bounded increasing $f, g : \{0,1\}^V \to \mathbb{R}$ we have*

$$\mu(fg) \geq \mu(f)\mu(g)\,.$$

(In parts of the probability literature, this is also known as *positive associations.*) Positive correlations holds for product measures – a result that goes back to Harris' classical paper [62] – while extremely useful criteria for positive correlations in dependent settings are provided by (various versions of) the FKG (Fortuin–Kasteleyn–Ginibre) inequality; see, e.g., [34]. The main result of Gandolfi et al. is the following.

Theorem 14.2 (Gandolfi, Keane and Russo [33]) *If $\mu$ is a translation invariant probability measure on $\{0,1\}^{\mathbb{Z}^2}$ satisfying (a) positive correlations, (b) invariance under interchange of or reflection in coordinate hyperplanes, and (c) ergodicity under horizontal translations and under vertical translations separately, then the number of infinite 1-clusters plus the number of infinite 0-clusters is $\mu$-a.s. at most 1.*

This generalizes the Ising model result of Cognilio et al. [26], and has later found interesting applications in – and adaptations to – the study of a wide variety of two-dimensional systems including the Potts model (Chayes [25], Häggström [49]), the hard-core model (Häggström [48]), the beach model (Hallberg [59]), and the random triangle model (Häggström and Jonasson [53]), plus an application to rigidity percolation that we will talk about in Section 16.

Here, instead of the somewhat involved geometric arguments needed to prove Theorem 14.2, we will settle for a shorter and remarkably elegant proof of the following variant of the result, where the condition of finite energy is added (which allows a slight weakening of the ergodicity assumptions). The finite energy condition holds in several (but not all) of the aforementioned applications. Other variants exist, including Sheffield's [97, Thm. 9.3.1] recent improvement where condition (b) of Theorem 14.3 is shown to be superfluous.

Theorem 14.3 *If $\mu$ is a translation invariant probability measure on $\{0,1\}^{\mathbb{Z}^2}$ satisfying (a) positive correlations, (b) invariance under interchange of or reflection in coordinate hyperplanes, (c') ergodicity, and (d) finite energy, then the number of infinite 1-clusters plus the number of infinite 0-clusters is $\mu$-a.s. at most 1.*

*Proof.* Following Georgii et al. [34], we will employ a beautiful geometric argument originally due to Yu Zhang [104]. First note that Theorem 12.2 adapts with zero additional complication to the setting where bond percolation is replaced by site percolation, so that we know that there is a.s. at most one infinite 1-cluster and at most one infinite 0-cluster.

Writing $A$ (resp. $B$) for the event that there is a unique infinite 1-cluster (resp. 0-cluster), we may thus assume, for contradiction, that $\mu(A \cap B) > 0$. By ergodicity, we then have

$$\mu(A \cap B) = 1 . \tag{14.1}$$

Take $\Lambda_n$ to denote $\{-n, \ldots, n\}^2$, i.e., $\Lambda_n$ is the box of side-length $2n+1$ centered at the origin. Defining $A_n$ (resp. $B_n$) as the event that the unique infinite 1-cluster (resp. 0-cluster) intersects $\Lambda_n$, we have from (14.1) that we can find an $n$ sufficiently large so that

$$\mu(A_n \cap B_n) \geq 0.999 .$$

Take $A_n^L$ (resp. $A_n^R, A_n^T, A_n^B$) to denote the event that for some vertex $x$ on the left (resp. right, top, bottom) side of the boundary of $\Lambda_n$ we can find an infinite self-avoiding path of 1's starting at $x$ but not containing any other vertex in $\Lambda_n$. Define $B_n^L, B_n^R, B_n^T$ and $B_n^B$ analogously, and write $\neg A_n^L$ for the complement of the event $A_n^L$, etc. Note that

$$A_n = A_n^L \cup A_n^R \cup A_n^T \cup A_n^B ,$$

and also that all four events on the right-hand side are increasing. Now we use the usual square root trick. By the assumed positive correlations property (a), we thus have

$$\begin{aligned} \mu(A_n) &= \mu(A_n^L \cup A_n^R \cup A_n^T \cup A_n^B) \\ &= 1 - \mu(\neg A_n^L \cap \neg A_n^R \cap \neg A_n^T \cap \neg A_n^B) \\ &\leq 1 - \mu(\neg A_n^L)\mu(\neg A_n^R)\mu(\neg A_n^T)\mu(\neg A_n^B) . \end{aligned}$$

By condition (b) all four events on the last line have the same $\mu$-probability, and we can deduce that

$$\mu(\neg A_n^L) \leq (1 - \mu(A_n))^{1/4}$$

and that

$$\begin{aligned} \mu(A_n^L) = \mu(A_n^R) &\geq 1 - (1 - \mu(A_n))^{1/4} \\ &\geq 1 - 0.001^{1/4} \\ &> 0.82 . \end{aligned} \tag{14.2}$$

By the same token, we get

$$\mu(B_n^T) = \mu(B_n^B) > 0.82 . \tag{14.3}$$

Finally, we define the event $D = A_n^L \cap A_n^R \cap B_n^T \cap B_n^B$, and use (14.2) and (14.3) to conclude that

$$\begin{aligned} \mu(D) &\geq 1 - 4(1 - 0.82) \\ &= 0.28 && (14.4) \\ &> 0 . && (14.5) \end{aligned}$$

Now consider the geometry of the infinite clusters on the event $D$. The path witnessing the event $A_n^L$ must of course be in an infinite 1-cluster, and likewise for the path witnessing $A_n^R$. But by uniqueness of the infinite 1-cluster these infinite clusters must be one and the same. Note, crucially, that this infinite 1-cluster blocks the two infinite 0-clusters witnessing $B_n^T$ and $B_n^B$ from merging, so that there must be at least two infinite 0-clusters. Hence, (14.5) contradicts the a.s. uniqueness of infinite 0-clusters, and we are done. □

## 15. Entanglement

Consider iid bond percolation on $\mathbb{Z}^3$. For $p \in (0, 1)$, we may imagine the occurrence of disjoint connected components that can nevertheless not be separated or "pulled apart". As a concrete example, consider two connected components each consisting of a circuit of eight vertices, the first one going through the vertices

$$(0,0,0),\ (0,0,1),\ (0,0,2),\ (0,1,2),\ (0,2,2),\ (0,2,1),\ (0,2,0),\ (0,1,0)$$

and back to $(0,0,0)$, while the second goes through

$$(-1,1,1),\ (0,1,1),\ (1,1,1),\ (1,1,2),\ (1,1,3),\ (0,1,3),\ (-1,1,3),\ (-1,1,2)$$

and back to $(-1,1,1)$. These two connected components form a single *entangled* component.

In *entanglement percolation* – studied first in the physics literature by Kantor and Hassold [73] and later in mathematics by Aizenman and Grimmett [2], Holroyd [65, 67], Grimmett and Holroyd [41] and Häggström [50] – the model is the usual iid bond percolation, but focus is shifted from connected components to entangled components. The phenomenon of entanglement is (like knot theory) intrinsically three-dimensional, and for this reason the study of entanglement percolation has so far dealt exclusively with the $\mathbb{Z}^3$ lattice.

Before we can proceed with the results, a series of definitions is needed. By a *sphere* in $\mathbb{R}^3$, we mean a subset of $\mathbb{R}^3$ homeomorphic to the unit sphere $\{x \in \mathbb{R}^3 : \|x\| = 1\}$, where $\|\cdot\|$ denotes Euclidean norm; for technical reasons, we also need to impose the condition that a sphere is a compact union of a finite number of polyhedral pieces. The complement of a sphere $S$ has precisely two connected components, called the *inside* and the *outside* (in the obvious way). For a sphere $S$ and an arbitrary set $A \subset \mathbb{R}^3$, $S$ is said to *separate* $A$ if $A$ intersects both the inside and the outside of $S$, but not $S$ itself.

We will identify any edge in $\mathbb{Z}^3$ with the unit length closed line segment connecting its endpoints. For a subset $F$ of the edges of $\mathbb{Z}^3$, write $[F]$ for the union of the corresponding line segments.

For $F$ finite, it is now clear what we ought to mean for $F$ to be entangled, namely that no sphere exists which separates $[F]$. Extending this to infinite $F$ can be done in several reasonable but non-equivalent ways; see [41]. For definiteness, we settle here for the following. An infinite edge set $F$ is said to be entangled if for every finite $F' \subset F$ there exists another finite set $F''$ such that $F' \subseteq F'' \subset F$ and $F''$ is entangled.

For iid bond percolation on $\mathbb{Z}^3$ there exists (by the same arguments as those in Section 1 for the critical value $p_c$) a critical value $p_e \in [0,1]$ such that

$$P_p(\text{an infinite entangled component exists}) = \begin{cases} 0 & \text{if } p < p_e \\ 1 & \text{if } p > p_e\,. \end{cases} \tag{15.1}$$

Obviously, since a connected subgraph of $\mathbb{Z}^3$ is automatically entangled, we have $p_e \leq p_c$. Aizenman and Grimmett [2] strengthened this to the strict inequality $p_e < p_c$, and later Holroyd [65] obtained the (surprisingly difficult, although see [29] for a striking example of a similar and "equally obvious" statement that turns out to be false) result that $p_e$ is strictly positive. Combining these results, we thus get

$$0 < p_e < p_c\,.$$

Numerical simulations (see [73]) suggest a rather narrow gap between $p_e$ and $p_c$:

$$p_c - p_e \approx 1.8 \cdot 10^{-7}\,. \tag{15.2}$$

As to the issue of uniqueness of the infinite entangled component, Grimmett and Holroyd [41] demonstrated that

$$P_p(\exists \text{ a unique infinite entangled component}) = 1 \tag{15.3}$$

for all $p > p_c$. Häggström [50] then improved this by showing that (15.3) holds for all $p > p_e$. In view of (15.2), this sounds like a rather unimpressive improvement, but it is in line with the general (and well-justified) obsession among percolation-theorists to prove properties of supercritical percolation not just for $p$ large but all the way down to the critical value.

A few words about the proof techniques employed for the uniqueness results in [41] and [50] are in order. The Burton–Keane trifurcation argument outlined in Section 2 is an extremely powerful idea, and would certainly be a very natural first attempt at proving (15.3). Difficulties arise, however, when trying to carry out the local modifications that form a central ingredient in the Burton–Keane approach. This is due to the highly erratic and non-local[7] character of entanglement (compared to connectivity). For instance, removing a single edge can cause an infinite entangled component to fall apart into infinitely many finite entangled components; see [65] or [41] for striking examples.

In order to prove (15.3) for $p > p_c$, Grimmett and Holroyd [41] instead chose to exploit a beautiful result of Barsky et al. [7] to the effect that, when $p > p_c(\mathbb{Z}^3)$, an infinite connected component arises even when restricting the percolation process to the upper half-space in $\mathbb{Z}^3$. A kind of sequential exploration (not entirely unlike the one in the proof of Theorem 5.6) of the entangled component containing a given vertex $x$ is defined, and every time this exploration process enters a hitherto untouched half-space, it has a fresh non-zero probability of hitting an infinite connected component of this half-space. It follows that if the entangled component containing $x$ is infinite, then it must a.s. intersect the infinite connected component in $\mathbb{Z}^3$ (which we know from Section 2 is unique), and therefore contain it. But obviously only one entangled component can do so, whence (15.3) follows.

This approach obviously does not work for $p \in (p_e, p_c)$ since then there is no infinite connected component to hit. Häggström [50] instead found an elaborate way to overcome the difficulties with adapting the Burton–Keane approach. A bird's-eye view of his proof is as follows. First, a uniqueness monotonicity result for infinite entangled components was obtained, borrowing ideas from [56] (see Section 5). It follows that if uniqueness fails somewhere on $(p_e, 1]$, it must do so on an interval of nonzero length. Assuming (for contradiction) that uniqueness fails, we can therefore find two different values $p_1 < p_2$ at which it fails. An easy adaptation of Lemma 2.6 shows that there must be infinitely many infinite entangled components at $p = p_1$. Using Coupling 1.1 and a suitable way of merging entangled components in $X_{p_1}$ depending on what happens in $X_{p_2}$ (defined as a kind of "interpolation between connectivity and entanglement", specifically designed to avoid the non-local behavior mentioned above) yields a configuration whose components exhibit analogs of trifurcations, and a contradiction can be deduced as in the original Burton–Keane argument. See [50] for the details.

[7]The entanglement concept considered here should of course not be confused with the one in quantum mechanics, despite the fact that the term "non-local" comes up in both contexts.

## 16. Rigidity

Yet another physically interesting property of percolation configurations, besides entanglement and connectivity, is *rigidity.* While entanglement is a weaker notion than connectivity (a connected graph is automatically entangled), rigidity is stronger. Imagine each edge in a graph as a hard rod; roughly speaking, the graph is rigid if it cannot be deformed via smooth movements of its vertices.

Rigidity is a dimension-dependent concept. Consider for instance four vertices in generic position in $\mathbb{R}^2$ or $\mathbb{R}^3$, connected to each other in a cycle. The resulting structure is non-rigid, because two non-adjacent vertices may be pulled apart or pushed towards each other. Adding a diagonal to the four-cycle makes the graph rigid in two dimensions, but not in three dimensions where it can still be bent (a bit like opening a book) along the diagonal edge.

For a rigorous definition, we proceed as follows (see, e.g., Holroyd [64] for more detail). Let $G = (V, E)$ be a finite graph, and fix $d \geq 2$. An *embedding* of $G$ in $\mathbb{R}^d$ is a map $r : V \to \mathbb{R}^d$, and the pair $(G, r)$ is called a *framework.* A $d$-dimensional *motion* of $(G, r)$ is a differentiable family $\{r_t : t \in [0, 1]\}$ of embeddings of $G$ in $\mathbb{R}^d$ such that for all $x, y \in V$ sharing an edge $\{x, y\} \in E$, and all $t$, we have

$$\|r_t(x) - r_t(y)\| = \|r(x) - r(y)\| . \tag{16.1}$$

The motion is *rigid* if (16.1) holds for all $x, y \in V$ (not just those pairs sharing an edge), and the framework $(G, r)$ is said to be ($d$-dimensionally) rigid if all its ($d$-dimensional) motions are rigid.

Whether the framework $(G, r)$ is rigid turns out to depend not only on $G$ but also on the embedding $r$. The good news, however, is that either almost all or almost no (with respect to $|V|^d$-dimensional Lebesgue measure) embeddings are rigid, and we say that $G$ is $d$-dimensionally rigid if almost all its embeddings in $\mathbb{R}^d$ are rigid.

As with entanglement, there are some options in how to extend the notion of rigidity to infinite graphs. Analogously to what we did in Section 15, we opt for the following. An infinite graph $G = (V, E)$ is said to be rigid if every finite subgraph of $G$ is contained in some rigid finite subgraph of $G$.

Rigidity in the context of percolation was first studied by Jacobs and Thorpe [68, 69]. In two dimensions, it is pointless from a rigidity point of view to work with the $\mathbb{Z}^2$ lattice, which is in itself non-rigid (and contains no rigid subgraphs with more than a single edge); more interesting is the triangular lattice $\mathbb{T}$, which is rigid. For iid bond percolation on $\mathbb{T}$, we have, analogously to (15.1), the existence of a rigidity critical value $p_r \in [0, 1]$ such that

$$P_p(\text{an infinite rigid component exists}) = \begin{cases} 0 & \text{if } p < p_r \\ 1 & \text{if } p > p_r\,. \end{cases}$$

Obviously, $p_c \le p_r \le 1$. Holroyd improved this in both ends, to

$$p_c < p_r < 1\,,$$

and in fact showed this for a large class of other lattices in two and higher dimensions; the main proviso is that the lattice itself is rigid. Unlike for entanglement (recall (15.2)), the gap $p_r - p_c$ appears (for $\mathbb{T}$) to be substantial: $p_c(\mathbb{T}) = 2\sin(\pi/18) \approx 0.347$ [103], while simulations in [69] suggest that $p_r(\mathbb{T}) \approx 0.660$.

Concerning uniqueness of the infinite rigid component, Holroyd [64] showed for iid percolation on $\mathbb{T}$ that

$$P_p(\exists \text{ a unique infinite rigid component}) \;=\; 1$$

for all $p \in (p_r, 1]$ *with at most countably many exceptional* $p$. This may to the inexperienced reader sound like a rather bizarre result, but in fact results with precisely this kind of proviso are not uncommon in statistical mechanics; see, e.g., [14, Thm 4.1] or [55, Thm 1.2]. And it does sound less strange with the following clarification: Write $\theta_r(p)$ for the $P_p$-probability that the origin is in an infinite rigid component. Holroyd established uniqueness of the infinite rigid component for all $p > p_r$ such that $\theta_r(p)$ is continuous at $p$. Note now that $\theta_r(p)$ is obviously increasing, and that an increasing function can have at most countably many discontinuities.

Holroyd's uniqueness result involved an adaptation of the Burton–Keane trifurcation approach. Similarly to the case of entanglement percolation in the previous section, the adaptation was a highly non-trivial task, stemming from the fact that rigidity, like entanglement, behaves in a highly non-local manner: local changes to a configuration can have global effects in terms of rigidity that are far more difficult to control than those of connectivity.

Using planarity ideas based on the Gandolfi–Keane–Russo [33] uniqueness result that we discussed in Section 14, Häggström [51] was able to remove the possible countable exceptional set in Holroyd's result, and establish uniqueness of the infinite rigid component for all $p$ (although the setting was still restricted to the triangular lattice $\mathbb{T}$).

In fact, both Holroyd's and Häggström's approaches involved specifically planar features. The question of how to establish uniqueness of the infinite rigid components in higher dimensions (or even for non-planar lattices in dimension $d = 2$) therefore remained open for a couple of years, until Häggström in his next paper [52] on rigidity found a way to prove the desired uniqueness for all $p > p_r(G)$ on any $d$-dimensional lattice $G$ (suitably defined) that is in itself rigid. The overall approach in [52] is similar to that for entanglement percolation in [50] outlined at the end of Section 15 – first obtain a uniqueness monotonicity result along the lines of Section 5, and then, assuming non-uniqueness happens for some $p$, look at percolation at two different such parameter values $p_1$ and $p_2$ simultaneously in order to construct an "intermediate" configuration to which a Burton–Keane-type argument can be used to establish a contradiction – but the overall implementation of this program turns out to be a bit of a rough ride.

One shortcoming, finally, of the general result in [52] is that it is silent on what happens *at* the critical value $p = p_r$, in contrast to the earlier result in [51] for $\mathbb{T}$ where having more than one infinite rigid component at criticality is ruled out. (For rigidity as well as entanglement, the nonexistence of an infinite component at criticality is far less confidently expected compared to the classical conjecture of nonexistence of infinite connected components at $p = p_c$ for $\mathbb{Z}^d$, mentioned in Section 1.)

**Acknowledgement.** We are grateful to Geraldine Bosco, Russ Lyons, Jeff Steif, Johan Tykesson and an anonymous referee, for helpful comments and corrections on earlier drafts of this paper.

## References

[1] Aizenman, M., Chayes, J.T., Chayes, L. and Newman, C.M. (1988) Discontinuity of the magnetization in one-dimensional $1/|x-y|^2$ Ising and Potts models, *J. Statist. Phys.* **50**, 1–40. MR0939480

[2] Aizenman, M. and Grimmett, G.R. (1991) Strict monotonicity for critical points in percolation and ferromagnetic models, *J. Statist. Phys.* **63**, 817–835. MR1116036

[3] Aizenman, M., Kesten, H. and Newman, C.M. (1987) Uniqueness of the infinite cluster and continuity of connectivity functions for short- and long-range precolation, *Comm. Math. Phys.* **111**, 505–532. MR0901151

[4] Alexander, K.S. (1995a) Simultaneous uniqueness of infinite clusters in stationary random labeled graphs, *Comm. Math. Phys.* **168**, 39–55. MR1324390

[5] Alexander, K.S. (1995b) Percolation and minimal spanning forests in infinite graphs, *Ann. Probab.* **23**, 87–104. MR1330762

[6] Babson, E. and Benjamini, I. (1999) Cut sets and normed cohomology with application to percolation, *Proc. Amer. Math. Soc.* **127**, 589–597. MR1622785

[7] Barsky, D.J., Grimmett, G.R. and Newman, C.M (1991) Percolation in half-spaces: equality of critical densities and continuity of the percolation probability, *Probab. Th. Rel. Fields* **90**, 111–148. MR1124831

[8] Benjamini, I., Kesten, H., Peres, Y. and Schramm, O. (2004) Geometry of the uniform spanning forest: transitions in dimensions 4, 8, 12,..., *Ann. Math.* **160**, 465–491. MR2123930
[9] Benjamini, I., Lyons, R., Peres, Y. and Schramm, O. (1999a) Critical percolation on any nonamenable graph has no infinite clusters, *Ann. Probab.* **27**, 1347–1356. MR1733151
[10] Benjamini, I., Lyons, R., Peres Y. and Schramm, O. (1999b) Group-invariant percolation on graphs, *Geom. Funct. Analysis* **9**, 29–66. MR1675890
[11] Benjamini, I., Lyons, R., Peres Y. and Schramm, O. (2001) Uniform spanning forests, *Ann. Probab.* **29**, 1–65. MR1825141
[12] Benjamini, I. and Schramm, O. (1996) Percolation beyond $\mathbb{Z}^d$, many questions and a few answers, *Electr. Comm. Probab.* **1**, 71–82. MR1423907
[13] Benjamini, I. and Schramm, O. (2001) Percolation in the hyperbolic plane, *J. Amer. Math. Soc.* **14**, 487–507. MR1815220
[14] van den Berg, J. and Steif, J.E. (1994) Percolation and the hard-core lattice gas model, *Stoch. Proc. Appl.* **49**, 179–197. MR1260188
[15] Bollobás, B. (1998) *Modern Graph Theory*, Springer, New York. MR1633290
[16] Borgs, C., Chayes, J.T., van der Hofstad, R., Slade, G., and Spencer, J. (2005) Random subgraphs of finite graphs: I. The scaling window under the triangle condition, *Random Structures Algorithms*, to appear. MR2155704
[17] Borgs, C., Chayes, J.T., Kesten, H., and Spencer, J. (2001) The birth of the infinite cluster: finite-size scaling in percolation. *Comm. Math. Phys.* **224**, 153–204. MR1868996
[18] Broadbent, S.R. and Hammersley, J.M. (1957) Percolation processes I: Crystals and mazes, *Proc. Cambridge. Phil. Soc.* **53**, 629–641. MR0091567
[19] Burton, R.M. and Keane, M.S. (1989) Density and uniqueness in percolation, *Comm. Math. Phys.* **121**, 501–505. MR0990777
[20] Burton, R.M. and Keane, M.S. (1991) Topological and metric properties of infinite clusters in stationary two-dimensional site percolation, *Israel J. Math.* **76**, 299–316. MR1177347
[21] Campanino, M. and Russo, L. (1985) An upper bound on the critical percolation probability for the three-dimensional cubic lattice, *Ann. Probab.* **13**, 478–491. MR0781418
[22] Chaboud, T. and Kenyon, C. (1996) Planar Cayley graphs with regular dual, *Internat. J. Algebra Comput.* **6**, 553–561. MR1419130
[23] Chayes, J.T., Chayes, L. and Newman, C.M. (1985) The stochastic geometry of invasion percolation, *Comm. Math. Phys.* **101**, 383–407. MR0815191
[24] Chayes, L. (1995) Aspects of the fractal percolation process, in *Fractal geometry and stochastics (Finsterbergen, 1994)*, pp 113–143, *Progr. Probab.* **37**, Birkhäuser, Basel. MR1391973
[25] Chayes, L. (1996) Percolation and ferromagnetism on $\mathbb{Z}^2$: the $q$-state Potts cases, *Stoch. Proc. Appl.* **65**, 209–216. MR1425356
[26] Coniglio, A., Nappi, C.R., Peruggi, F. and Russo, L. (1976) Percolation

and phase transitions in the Ising model, *Comm. Math. Phys.* **51**, 315–323. MR0426745

[27] Durrett, R. (1991) *Probability: Theory and Examples*, Wadsworth & Brooks/Cole, Pacific Grove. MR1068527

[28] Edwards, R.G. and Sokal, A.D. (1988) Generalization of the Fortuin–Kasteleyn–Swendsen–Wang representation and Monte Carlo algorithm, *Phys. Rev.* D **38**, 2009–2012. MR0965465

[29] van Enter, A.C.D. (1987) Proof of Straley's argument for bootstrap percolation, *J. Statist. Phys.* **48**, 943–945. MR0914911

[30] Erdős, P. and Rényi, A. (1959) On random graphs, *Publicationes Mathematicae Debrecen* **6**, 290–297. MR0120167

[31] Fortuin, C.M. and Kasteleyn, P.W. (1972) On the random-cluster model. I. Introduction and relation to other models, *Physica* **57**, 536–564. MR0359655

[32] Gandolfi, A., Keane, M. and Newman, C.M. (1992) Uniqueness of the infinite component in a random graph with applications to percolation and spin glasses, *Probab. Th. Rel. Fields* **92**, 511–527. MR1169017

[33] Gandolfi, A., Keane, M. and Russo, L. (1988) On the uniqueness of the infinite occupied cluster in dependent two-dimensional site percolation *Ann. Probab.* **16**, 1147–1157. MR0942759

[34] Georgii, H.-O., Häggström, O. and Maes, C. (2001) The random geometry of equilibrium phases, *Phase Transitions and Critical Phenomena, Volume 18* (C. Domb and J.L. Lebowitz, eds), pp 1-142, Academic Press, London. MR2014386

[35] Grimmett, G.R. (1994) Percolative problems, *Probability and Phase Transition* (ed. G.R. Grimmett), Kluwer, Dordrecht, 69–86. MR1283176

[36] Grimmett, G.R. (1995) The stochastic random-cluster process, and the uniqueness of random-cluster measures, *Ann. Probab.* **23**, 1461–1510. MR1379156

[37] Grimmett, G.R. (1999) *Percolation*, Springer, Berlin. MR1707339

[38] Grimmett, G.R. (2003) The random-cluster model, in *Probability on Discrete Structures*, vol **110** of *Encyclopedia of Mathematical Sciences*, pp 73–123, Springer, Berlin. MR2023651

[39] Grimmett, G.R. (2005) Uniqueness and multiplicity of infinite clusters, in *Dynamics and Stochastics: Festschrift in Honor of Michael Keane*, IMS Lecture Notes-Monograph Series, pp 24–36.

[40] Grimmett, G.R. (2006) *The Random-Cluster Model*, Springer, Berlin. MR2243761

[41] Grimmett, G.R. and Holroyd, A.E. (2000) Entanglement in percolation, *Proc. London Math. Soc.* **81**, 485–512. MR1770617

[42] Grimmett, G.R. and Newman, C.M. (1990) Percolation in $\infty + 1$ dimensions, In *Disorder in Physical Systems* (G.R. Grimmett and D.J.A. Welsh, eds), pp. 167–190, Oxford University Press, New York. MR1064560

[43] Grossman, J.W. (2002) The evolution of the mathematical research collaboration graph, *Congressus Numerantium* **158**, 201–212. MR1985159

[44] Grossman, J.W. (2003) *The Erdős Number Project*,

http://www.oakland.edu/enp/

[45] Häggström, O. (1995) Random-cluster measures and uniform spanning trees, *Stoch. Proc. Appl.* **59**, 267–275. MR1357655

[46] Häggström, O. (1996) The random-cluster model on a homogeneous tree, *Probab. Th. Rel. Fields* **104**, 231–253. MR1373377

[47] Häggström, O. (1997a) Infinite clusters in dependent automorphism invariant percolation on trees, *Ann. Probab.* **25**, 1423–1436. MR1457624

[48] Häggström, O. (1997b) Ergodicity of the hard-core model on $\mathbf{Z}^2$ with parity-dependent activities, *Ark. Mat.* **35**, 171–184. MR1443040

[49] Häggström, O. (1999) Positive correlations in the fuzzy Potts model, *Ann. Appl. Probab.* **9**, 1149–1159. MR1728557

[50] Häggström, O. (2001a) Uniqueness of the infinite entangled component in three-dimensional bond percolation, *Ann. Probab.* **29**, 127–136. MR1825145

[51] Häggström, O. (2001b) Uniqueness in two-dimensional rigidity percolation, *Math. Proc. Cambridge Phil. Soc.* **130**, 175–188. MR1797779

[52] Häggström, O. (2003) Uniqueness of infinite rigid components in percolation models: the case of nonplanar lattices, *Probab. Th. Rel. Fields* **127**, 513–534. MR2021194

[53] Häggström, O. and Jonasson, J. (1999) Phase transition in the random triangle model, *J. Appl. Probab.* **36**, 1101–1115. MR1742153

[54] Häggström, O., Jonasson, J. and Lyons, R. (2002) Explicit isoperimetric constants and phase transitions in the random-cluster model, *Ann. Probab.* **30**, 443–473. MR1894115

[55] Häggström, O. and Pemantle, R. (2000) Absence of mutual unbounded growth for almost all parameter values in the two-type Richardson model, *Stoch. Proc. Appl.* **90**, 207–222. MR1794536

[56] Häggström, O. and Peres, Y. (1999) Monotonicity of uniqueness for percolation on transitive graphs: all infinite clusters are born simultaneously, *Probab. Th. Rel. Fields* **113**, 273–285. MR1676835

[57] Häggström, O., Peres, Y., and Schonmann, R.H. (1999) Percolation on transitive graphs as a coalescent process: relentless merging followed by simultaneous uniqueness, in *Perplexing Probability Problems: Papers in Honor of Harry Kesten* (M. Bramson and R. Durrett, eds), pp. 53–67, Birkhäuser, Boston. MR1703125

[58] Häggström, O., Peres, Y., and Steif, J. (1997) Dynamical percolation, *Ann. Inst. H. Poincaré, Probab. Stat.* **33**, 497–528. MR1465800

[59] Hallberg, P. (2004) *Gibbs Measures and Phase Transitions in Potts and Beach Models*, Ph.D. thesis, Royal Institute of Technology, Stockholm, http://media.lib.kth.se/dissengrefhit.asp?dissnr=3837

[60] Hammersley, J.M. (1961) Comparison of atom and bond percolation processes, *J. Math. Phys.* **2**, 728–733. MR0130722

[61] Hara, T. and Slade, G. (1994) Mean-field behaviour and the lace expansion, in *Probability and Phase Transition* (G. R. Grimmett, ed.), pp. 87–122, Kluwer Acad. Publ., Dordrecht. MR1283177

[62] Harris, T.E. (1960) A lower bound for the critical probability in a certain

percolation process, *Proc. Cambridge Phil. Soc.* **56**, 13–20. MR0115221
[63] Higuchi, Y. (1982) Coexistence of the infinite ∗-clusters: a remark on the square lattice site percolation, *Z. Wahrsch. Verw. Gebiete* **61**, 75–81. MR0671244
[64] Holroyd, A.E. (1998) Existence and uniqueness of infinite components in generic rigidity percolation, *Ann. Appl. Probab.* **8**, 944–973. MR1627815
[65] Holroyd, A.E. (2000) Existence of a phase transition for entanglement percolation, *Math. Proc. Cambridge Phil. Soc.* **129**, 231–251. MR1765912
[66] Holroyd, A.E. (2001) Rigidity percolation and boundary conditions, *Ann. Appl. Probab.* **11**, 1063–1078. MR1878290
[67] Holroyd, A.E. (2002) Entanglement and rigidity in percolation models, *In and Out of Equilibrium (Mambucaba, 2000)*, pp 299–307, *Progr. Probab.* **51**, Birkhäuser, Boston. MR1901959
[68] Jacobs, D.J. and Thorpe, M.F. (1995) Generic rigidity percolation: the pebble game, *Phys. Rev. Lett* **75**, 4051–4054.
[69] Jacobs, D.J. and Thorpe, M.F. (1996) Generic rigidity percolation in two dimensions, *Phys. Rev. E* **53**, 3682–2693.
[70] Janson, S., Łuczak, T. and Rucinski, A. (2000) *Random Graphs*, Wiley, New York. MR1782847
[71] Jonasson, J. (1999) The random cluster model on a general graph and a phase transition characterization of nonamenability, *Stoch. Proc. Appl.* **79**, 335–354. MR1671859
[72] Jonasson, J. and Steif, J. (1999) Amenability and phase transition in the Ising model, *J. Theor. Prob.* **12**, 549–559. MR1684757
[73] Kantor, T. and Hassold, G.N. (1988) Topological entanglements in the percolation problems, *Phys. Rev. Lett.* **60**, 1457–1460. MR0935098
[74] Kesten, H. (1959) Full Banach mean values on countable groups, *Math. Scand.* **7**, 146–156. MR0112053
[75] Kesten, H. (1959) Symmetric random walks on groups, *Trans. Amer. Math. Soc.* **92**, 336–354. MR0109367
[76] Kesten, H. (1980) The critical probability of bond percolation on the square lattice equals $\frac{1}{2}$, *Comm. Math. Phys.* **74**, 41–59. MR0575895
[77] Kesten, H. (1982) *Percolation Theory for Mathematicians*, Birkhäuser, Boston. MR0692943
[78] Lalley, S. (1998) Percolation on Fuchsian groups, *Ann. Inst. H. Poincaré, Probab. Stat.* **34**, 151–178. MR1614583
[79] Lyons, R (2000) Phase transitions on nonamenable graphs, *J. Math. Phys.* **41**, 1099–1126. MR1757952
[80] Lyons, R. and Peres, Y. (2005) *Probability on Trees and Networks*, Cambridge University Press, to appear, http://mypage.iu.edu/˜rdlyons/prbtree/prbtree.html
[81] Lyons, R., Peres, Y. and Schramm, O. (2006) Minimal spanning forests, *Ann. Probab.*, to appear.
[82] Lyons, R. and Schramm, O. (1999) Indistinguishability of percolation clusters, *Ann. Probab.* **27**, 1809–1836. MR1742889
[83] Meester, R. (1994) Uniqueness in percolation theory, *Statist. Neerl.* **48**,

237–252. MR1310339

[84] Meester, R. and Roy, R. (1996) *Continuum Percolation*, Cambridge University Press. MR1409145

[85] Mohar, B. (1988) Isoperimetric inequalities, growth and the spectrum of graphs, *Lin. Alg. Appl.* **103**, 119–131. MR0943998

[86] Newman, C.M. and Schulman, L.S. (1981) Infinite clusters in percolation models, *J. Statist. Phys.* **26**, 613–628. MR0648202

[87] Newman, C.M. and Stein, D.L. (1996) Ground-state structure in a highly disordered spin-glass model, *J. Statist. Phys.* **82**, 1113–1132. MR1372437

[88] Pak, I. and Smirnova-Nagnibeda, T. (2000) Uniqueness of percolation on nonamenable Cayley graphs, *Comptes Rendus Acad. Sci. Paris, Ser. I Math.* **330**, 495–500. MR1756965

[89] Pemantle, R. (1991) Choosing a spanning tree for the integer lattice uniformly, *Ann. Probab.* **19**, 1559–1574. MR1127715

[90] Peres, Y. (2000) Percolation on nonamenable products at the uniqueness threshold, *Ann. Inst. H. Poincaré, Probab. Stat.* **36**, 395–406. MR1770624

[91] Peres, Y., Pete, G. and Scolnicov, A. (2006) Critical percolation on certain nonunimodular graphs, *New York J. Math.* **12**, 1–18. MR2217160

[92] Peres, Y. and Steif, J.E. (1998) The number of infinite clusters in dynamical percolation, *Probab. Th. Rel. Fields* **111**, 141–165. MR1626782

[93] Pfister, C.-E. and Vande Velde, K. (1995) Almost sure quasilocality in the random cluster model, *J. Statist. Phys.* **79**, 765–774. MR1327908

[94] Schonmann, R.H. (1999) Percolation in $\infty+1$ dimensions at the uniqueness threshold, in *Perplexing Probability Problems: Papers in Honor of Harry Kesten* (M. Bramson and R. Durrett, ed.), pp. 53–67, Birkhäuser, Boston. MR1703124

[95] Schonmann, R.H. (1999) Stability of infinite clusters in supercritical percolation, *Probab. Th. Rel. Fields* **113**, 287–300. MR1676831

[96] Schramm, O. and Steif, J.E. (2005) Quantitative noise sensitivity and exceptional times for percolation, preprint.

[97] Sheffield, S. (2005) Random surfaces, *Astérisque* **304**, vi+175 pp. MR2251117

[98] Stacey, A.M. (1996) The existence of an intermediate phase for the contact process on trees, *Ann. Probab.* **24**, 1711–1726. MR1415226

[99] Timár, A. (2006a) Cutsets in infinite graphs, *Comb. Probab. Computing* **16**, 1–8.

[100] Timár, A. (2006b) Percolation on nonunimodular graphs, *Ann. Probab.*, to appear.

[101] Timár, A. (2006c) Neighboring clusters in Bernoulli percolation,*Ann. Probab.*, to appear.

[102] Trofimov, V.I. (1985) Automorphism groups of graphs as topological groups, *Math. Notes* **38**, 717–720. MR0811571

[103] Wierman, J.C. (1981) Bond percolation on honeycomb and triangular lattices, *Avd. Appl. Probab.* **13**, 298–313. MR0612205

[104] Zhang, Y. (1988) Unpublished, although see Grimmett [37, pp. 289–291].